\edef \helpstring{loaded}%
\ifx \headloaded\helpstring%
\else%
\magnification =\magstep0%
\font\micro =cmr12 at 8.00truept%
\font\SMALL =cmr12 at 9.00truept%
 at 10.00truept%
\font\small =cmr12 at 10.95truept%
 at 20.74truept
\font\large =cmr12%
\font\Large =cmr12 at 14.4truept%
\font\largeSC=cmcsc12%
\parindent=0pt%
\vsize = 21truecm%
\hsize = 14.2truecm%

\topskip =15 pt%
\voffset =1.5truecm%
\hoffset =0.7cm%
\rightskip = .0truecm%
\leftskip = .0truecm%
\footline ={\hss }%
\edef \headloaded{loaded}%
%
%
%
\edef\ONFLAG{on}%
\edef\OFFFLAG{off}%
\edef\YESFLAG{yes}%
\edef\NOFLAG{no}%
%
%
\edef\empty{}%
%
%
%
\newcount\VARX%
\newcount\VARY%
\newcount\RESULT%
\newcount\paragraphcounter
\paragraphcounter 0\relax%
%
%
%
%
%
%
\def\comparE #1#2{
    TT\fi%
    \edef\firsT{#1}%
    \edef\seconD{#2}%
    \ifx\firsT\seconD%
}%
%
%
%
%
%

%
%
%
%
%
%
%
%
%
\def\commentline #1{%
\smallskip\hbox to \hsize{%
\vrule height .4pt width 5cm depth 0.0 pt\hss%
{\sc #1}%
\hss\vrule height .4pt width 5cm depth 0.0 pt}}%

\long\def\comment #1{}
\def\endcomment{}
\long\def\begincomment #1{}
%
\newcount\longversiondepth%
\newcount\longversioncommentdepth%
\longversiondepth 0\relax%
\longversioncommentdepth 0\relax%
\edef\longversionflag{\OFFFLAG}%
\edef\longversioncommentflag{\OFFFLAG}%
\def\longversion #1#2{
  \comparedepth {\longversiondepth}{#1}%
  \ifx\result\NOFLAG%
  \else%
     \edef\longversionflag{\result}%
  \fi%
  \comparedepth {\longversioncommentdepth}{#2}%
  \ifx\result\NOFLAG%
  \else%
    \edef\longversioncommentflag{\result}%
  \fi%
}%
\long\def\beginlong #1{
%
  \ifx\longversionflag\OFFFLAG\relax%
  \else%
    \ifx\longversioncommentflag\OFFFLAG%
    \else%
      \par\zl{================== LONG VERSION ===================}\par%
    \fi%
    #1%
    \ifx\longversioncommentflag\OFFFLAG%
    \else%
      \par\zl{=============== END OF LONG VERSION ===============}\par%
    \fi%
  \fi%
}%
\def\endlong{}%
\def\longx#1{
	\ifx\longversionflag\OFFFLAG\relax%
	\else%
	    \ifx\longversioncommentflag\OFFFLAG%
    	\else%
	      @%
    	\fi%
	    #1%
	    \ifx\longversioncommentflag\OFFFLAG%
    	\else%
	      @%
    	\fi%
	\fi%
}%
\long\def\lversion #1{
  \ifx\longversionflag\OFFFLAG\relax%
  \else%
    \ifx\longversioncommentflag\OFFFLAG%
    \else%
      \par\zl{================== LONG VERSION ===================}\par%
    \fi%
    #1%
    \ifx\longversioncommentflag\OFFFLAG%
    \else%
      \par\zl{=============== END OF LONG VERSION ===============}\par%
    \fi%
  \fi%
}%
%
%
\def\text#1{{\leavevmode\hbox{#1}}}%
\def \newpage{\vfill\eject}%
\def \Paragraph{\S }%
\def \hcenter #1#2{
\setbox\helpbox =\hbox{#1}%
\setbox\helpboxb =\hbox{#2}%
\ifnum \wd\helpboxb > \wd\helpbox \box\helpboxb%
\else\hbox to \wd\helpbox{\hss #2\hss}%
\fi%
}%
\long\def\abstract #1{%
 \zl{\vbox {%
           \hsize = 0.9\hsize%
           {\small%
            {\smallsc Abstract}.%
            #1%
           }%
 }}%
 \bigskip%
}%
\def\OP #1#2{
{\bf Open problems.}\\ %
\beg {#1}%
#2%
\endbeg\medskip%
}%

%
%
%
%
\chardef\bodY\catcode`@\catcode`@=11%
\def\hexnumber@#1{\ifcase#1 0\or 1\or 2\or 3\or 4\or 5\or 6\or 7\or 8\or%
 9\or A\or B\or C\or D\or E\or F\fi}%
\edef\newfam{}
\def\newfam{\alloc@8\fam\chardef\sixt@@n}
\def \newfamilY #1#2{
  \if\comparE{#1}{}%
    \edef\lokaL{#2}%
  \else%
    \edef\lokaL{#1}%
  \fi%
  \expandafter\font\csname ten#2\endcsname =\lokaL 10%
  \expandafter\font\csname seven#2\endcsname =\lokaL 7%
  \expandafter\font\csname five#2\endcsname =\lokaL 5%
  \expandafter\newfam\csname #2fam\endcsname%
  \expandafter\textfont\csname #2fam\endcsname =\csname ten#2\endcsname%
  \expandafter\scriptfont\csname #2fam\endcsname =\csname seven#2\endcsname%
  \expandafter\scriptscriptfont\csname #2fam\endcsname =\csname five#2\endcsname%
  \expandafter\edef\csname #2fam@\endcsname{%
  \expandafter\hexnumber@\csname #2fam\endcsname}%
}%
\def\largeletteralternativE #1{\ifcase #1{A}{41}\or {B}{42}\or {C}{43}%
\or {D}{44}\or {E}{45}\or {F}{46}\or {G}{47}\or {H}{48}\or {I}{49}%
\or {J}{4A}\or {K}{4B}\or {L}{4C}\or {M}{4D}\or {N}{4E}\or {O}{4F}%
\or {P}{50}\or {Q}{51}\or {R}{52}\or {S}{53}\or {T}{54}\or {U}{55}%
\or {V}{56}\or {W}{57}\or {X}{58}\or {Y}{59}\or {Z}{5A}\fi}%
\def\letterlargE#1{\ifcase #1\or A\or B\or C%
\or D\or E\or F\or G\or H\or I%
\or J\or K\or L\or M\or N\or O%
\or P\or Q\or R\or S\or T\or U%
\or V\or W\or X\or Y\or Z\fi}%
\def\smallletteralternativE #1{\ifcase #1{a}{61}\or {b}{62}\or {c}{63}%
\or {d}{64}\or {e}{65}\or {f}{66}\or {g}{67}\or {h}{68}\or {i}{69}%
\or {j}{6A}\or {k}{6B}\or {l}{6C}\or {m}{6D}\or {n}{6E}\or {o}{6F}%
\or {p}{70}\or {q}{71}\or {r}{72}\or {s}{73}\or {t}{74}\or {u}{75}%
\or {v}{76}\or {w}{77}\or {x}{78}\or {y}{79}\or {z}{7A}\fi}%
\def\lettersmalL#1{\ifcase #1\or a\or b\or c%
\or d\or e\or f\or g\or h\or i%
\or j\or k\or l\or m\or n\or o%
\or p\or q\or r\or s\or t\or u%
\or v\or w\or x\or y\or z\fi}%
\def\readfirstentrY#1#2{#1}%
\def\readsecondentrY#1#2{#2}%
\def\mathcharacterS#1#2#3{
  \VARX 0\relax%
  \loop%
    \edef\lokaL{\csname #3letteralternativE\endcsname {\VARX}}\relax%
    \edef\executE{\csname #1\expandafter\readfirstentrY\lokaL\endcsname%
      "0\csname #2fam@\endcsname\expandafter\readsecondentrY\lokaL}\relax%
    \expandafter\mathchardef\executE%
  \ifnum\VARX <25%
    \advance\VARX +1\relax%
  \repeat%
}%
%
%
\newfamilY{msam}{msa}%
\newfamilY{msbm}{msb}%
\newfamilY{}{rsfs}%
\newfamilY{}{eufm}%
\newfamilY{}{eusm}%
\newfamilY{}{cmcsc}%
\edef\kalifam@{2}%
\mathcharacterS{L}{msb}{large}%
\mathcharacterS{S}{rsfs}{large}%
\mathcharacterS{D}{eufm}{large}%
\mathcharacterS{D}{eufm}{small}%
\mathcharacterS{C}{kali}{large}%
%
%
\def \sc{\fam \z@ \tencmcsc}
\font \smallsc =cmcsc10 at 10.00truept%

%
%
%
%
%
%
%
\mathchardef\subseteqq"3\msafam@6A%
\mathchardef\supseteqq"3\msafam@6B%
\mathchardef\nmid"3\msbfam@2D%
\mathchardef\proper"3\msbfam@24%
\mathchardef\prober"3\msbfam@25%
\mathchardef\square"3\msafam@03%
\mathchardef\blacksquare"3\msafam@04%
\mathchardef\bigstar"3\msafam@46%
\mathchardef\constant"3\eusmfam@43%
\mathchardef\ltr"3\msafam@43%
\mathchardef\rtr"3\msafam@42%
%
%
%
%
%
\mathchardef\uhr"3\msafam@16%
\mathchardef\spez"3\msafam@20%
\mathchardef\semidirect"3\msbfam@6E%
\mathchardef\ltimes"3\msbfam@6E%
\mathchardef\semiDirect"3\msbfam@6F%
\mathchardef\rtimes"3\msbfam@6F%
\mathchardef\sqleq"3276%
\mathchardef\sqgeq"3277%
\mathchardef\upharpoonright"3225%
\mathchardef\downharpoonright"3226%
\mathchardef\heart"327E%
\def \herz #1{($\heart%
  \if\comparE{#1}{}%
    leer
  \else%
    \VARX #1\relax
    \ifcase \VARX%
    \or%
    \or%
      \heart%
    \or%
      \heart\heart%
    \fi%
  \fi%
  $)\relax%
}%
%
%
%
%
%
%
%
\catcode`@\bodY%
%
%
%
%
\def\R {{\rm{I\kern -0.17 em R}}}%
\def\N {{\rm{I\kern -0.17 em N}}}%
\def\F {{\rm{I\kern -0.17 em F}}}%
\def\makeQC#1#2#3#4#5#6{%
  {\rm{ \kern #1\vrule height #2 depth #3 width #4 \kern #5 #6}}%
}%
\edef\Q{{\mathchoice
  {\makeQC{0.26 em}{1.6 ex}{-0.1 ex}{0.05 em}{-0.31 em}{Q}}%
  {\makeQC{0.26 em}{1.6 ex}{-0.1 ex}{0.05 em}{-0.31 em}{Q}}%
  {\makeQC{0.24 em}{1.1 ex}{-0.05 ex}{0.05 em}{-0.29 em}{Q}}%
  {\makeQC{0.22 em}{.75 ex}{-0.05 ex}{0.05 em}{-0.27 em}{Q}}%
}}%

%
%
%
%
%
\newcount\factnumber
\newcount\microfactnumber%
\newcount\exercisenumber%
\newcount\pointnumber%
\newcount \listlevel
\def \setcount#1#2{#1#2\relax}
\def\stepcount#1{\advance#1 1\relax}
%
%
%
%
\newbox \listentry
\newbox\helpbox
\newbox\helpboxb%
%
%
%
%
\def \\ {%
        \ifmmode%
                \bar%
        \else%
                \par%
        \fi}%
\def\einzug {\\ \hglue 3mm}%
%
%
%
%
\def \chaptertitle #1#2#3{%
\ifdim \the\pagetotal > 0pt\newpage%
\else\fi%
\global\def \chapternumber {#1}%
\global\def \chaptername {#3}%
{\Large Chapter #1: #2}\bigskip%
}%
\def \midheadlineleft{{\sc Marcus Tressl}}
%
\def \doubleflag {double}
\def \singleflag {single}%
\newcount\headlineshapedepth%
\headlineshapedepth 0\relax%
\edef\headlineshape{\OFFFLAG}%
\def \headlineshape #1#2#3{
%
%
%
  \comparedepth{\headlineshapedepth}{#1}%
  \ifx\result\ONFLAG%
    \def \titleboxflag {#2}%
    \def \layoutflag {#3}%
  \else%
  \fi%
}%
\def \titlebox #1#2{
 \if\comparE{#1}{\OFFFLAG}%
    \if\comparE{#2}{}%
     \else%
       \Paragraph #2%
     \fi%
 \else%
   {\micro \vbox {%
                \hcenter{\the\day .\the\month .\the\year}{%
                         \if\comparE{#2}{}%
                         \else\Paragraph #2%
                         \fi%
                         {\SMALL \jobname }%
                }\relax%
                \vskip -2ex\hbox{\the\day .\the\month .\the\year}%
         }\relax%
   }%
 \fi%
}%
\def\myheadline#1#2#3{
    \vbox{\hbox to \hsize{{\rm #1\hss #2\hss #3}}\relax%
          }\relax
}%
\def\headlinetextvalue{}%
\def \constructheadline #1#2#3#4#5{
  \if\comparE{#5}{}%
    \def \midheadline {\headlinetextvalue}%
  \else%
    \def \midheadline {#5}%
  \fi%
  \ifx #2\singleflag%
    \headline={\relax%
      \myheadline {\titlebox{#1}{#3}}{\midheadline}{\folio}%
    }%
  \else\fi%
  \ifx #2\doubleflag%
    \headline={\relax%
      \ifodd\pageno%
        \ifnum\pageno=1%
          \hss%
        \else%
          \myheadline {\titlebox{#1}{#3}}{\midheadline}{\folio}%
        \fi%
      \else%
        \myheadline {\folio}{\midheadlineleft}{\titlebox{#1}{#3}}%
      \fi%
    }%
  \else\fi%
}%
\def\titlewire #1{
 \ifdim \the\pagetotal > 650pt \newpage
 \else  \ifdim \the\pagetotal < 10pt
        \else \vglue 1cm
        \fi%
 \fi%
 \setcount\factnumber 0%
#1\bigskip\medskip%
}%
\def \title #1#2#3#4{
  \def \parmone {on}
  \def \parmtwo {single}
  \edef \lokaL {#1}
  \ifx \lokaL\OFFFLAG%
       \def \parmone {off}\else\fi%
  \ifx \lokaL\doubleflag%
       \def \parmone {off}%
       \def \parmtwo {double}\else\fi%
  \ifx \titleboxflag\OFFFLAG%
       \def \parmone {off}\else\fi%
  \ifx \titleboxflag\ONFLAG%
       \def \parmone {on}\else\fi%
  \ifx \layoutflag\doubleflag%
       \def \parmtwo {double}\else\fi%
  \ifx \layoutflag\singleflag%
       \def \parmtwo {single}\else\fi%
  \constructheadline{\parmone}{\parmtwo}{}{#3}{#4}%
   \if\comparE{#2}{}%
    \global\def \paragraphnumber{}%
    \titlewire{{\bf \large #3}}%
  \else%
  	\if\comparE{#2}{*}%
  		\advance \paragraphcounter 1\relax%
  		\global\def \paragraphnumber{\the \paragraphcounter .}%
    	\titlewire{\zl{\largeSC \paragraphnumber \ #3}}%
  	\else%
    	\global\def \paragraphnumber{#2.}%
    	\titlewire{\zl{\largeSC  #2. #3}%
		}%
	\fi%
  \fi%
}%
%
%
\def \subtitle #1{%
  \ifdim \the\pagetotal > 500pt \newpage
  \else  \ifdim \the\pagetotal < 10pt
        \else \vglue .7cm
        \fi%
  \fi%
  \zl {{\bf #1}}\bigskip
}%
%
%
%
%
%
%
%
%
\def \item#1{%
\advance \listlevel by 1\relax%
\par%
\noindent%
\hbox to 0pt{%
\hskip -\wd\listentry%
#1%
\hfil%
}}%
\def \ro {
  \item {(\romannumeral\listlevel )}%
}%
\def \ar {
  \item {(\number\listlevel)}%
}%
\def \lsmall {
  \item {(\lettersmalL{\listlevel})}%
}%
\def \beg #1{%
\par%
\bgroup%
\listlevel=0%
\setbox\listentry=\hbox{#1\hskip 0.3em}
\advance \leftskip by \wd\listentry \relax%
}%
\def \endbeg{%
\par%
\egroup%
}%
%
%
%
%
\def\Punkt {\bigskip%
            \stepcount{\pointnumber}%
            {\bf \the \pointnumber .}}%
\def\punkt {\the\pointnumber .}%
%
%
\def\factspacevalue{}%
\long\def \fact#1#2#3{\edef\factname{#1 }%
\factspacevalue
\stepcount{\factnumber}%
 (\paragraphnumber \the\factnumber )\ {\sc #1.}\label{#2}{%
 (\paragraphnumber \the\factnumber )}%
 \ifx \showlabelflag\compareshowlabel \\ %
   {(\paragraphnumber \the\factnumber )\ \sc #1.}%
 \else\fi%
 {\it #3}%
}%
\def \ende {\hglue 0pt$\hfill \square $\medskip\par \penalty -200}%
\def \underconstruction #1{\\ {\bf UNDER CONSTRUCTION:} {\it {#1}}\relax \ende}%
\def \pr {\medskip {\sc Proof. }}%
\long\def \Definition #1{\stepcount{\factnumber}%
(\paragraphnumber \the\factnumber )\ {\sc Definition.}\label{#1}{%
(\paragraphnumber \the\factnumber )}}%
\long\def \Remark #1{\stepcount{\factnumber}%
(\paragraphnumber \the\factnumber )\ {\sc Remark.}\label{#1}{%
\paragraphnumber \the\factnumber }}%
\def \claim {\ub {Claim.} }%
\def \Claim #1{\ub {Claim #1.} }%
\def \case #1{\ub {Case #1.} }%
%
%
%
%
%
%
%
%
%
%
%
%
%
%
%
%
\def \dcl{\mathop {\rm dcl}\nolimits}%
%
%
%
%
%
%
%
%
%
%
%
%
%
%
%
%
%
%
%
%
%
%
%
%
%
%
%
%
%
%
%
%
%
%
%
%
%
%
%
%
%
%
%
%
%
%
%
\def \sign {\mathop {\rm sign}\nolimits }%
\def \Wl {{\raise2.3pt\hbox{%
${\scriptscriptstyle \not \ }$}}\hskip -0.6ex{\rm L}}%
\def \precc#1{\! \! \mathrel{\mathop{\kern 0pt \prec}\limits _{\ \ \ #1}}}%
%
%
%
%
\def \mal{\! \cdot \! }%
\def \valO {{O\hskip -2mm^\smile }}%
%
%
%
%
%
%
\def\chix{{\raise 2pt\hbox{$\chi $}}}%
%
%
%
\def \<{\langle}%
\def \>{\rangle}%
\def \phi {\varphi }%
\def \theta {\vartheta }%
\def \dnf {\mathop {{\raise -3.9pt\hbox{$\smile $}} \hskip -1.5ex \vert%
\hskip 1ex}}%
\def \eps{\varepsilon}%
\def \lra{\longrightarrow}%
\def \ub #1{{\underbar {#1}}}%
\def \zl #1{\centerline {#1}}%
\def \a {\vert}%
\def \iff {{\ \Leftrightarrow \ }}%
\def \st {\ \vert \ }%
\def \notion #1{{\bf #1}}%
%
%
%
%
\def\stapel #1#2#3#4#5{
\setbox \helpbox=\hbox{$#1$}%
\setbox \helpboxb=\hbox{$#2$}%
\ifnum \wd\helpboxb > \wd\helpbox \setbox \helpbox=\hbox{$#2$}\else\fi%
\setbox \helpboxb=\hbox{$#3$}%
\ifnum \wd\helpboxb > \wd\helpbox \setbox \helpbox=\hbox{$#3$}\else\fi%
\setbox \helpboxb=\hbox{$#4$}%
\ifnum \wd\helpboxb > \wd\helpbox \setbox \helpbox=\hbox{$#4$}\else\fi%
\vcenter {%
\def\helpmacro{#1}%
\ifx \helpmacro\empty%
\else \hcenter{\box\helpbox}{$#1$}\fi%
\def\helpmacro{#2}%
\ifx \helpmacro\empty%
\else\vskip -#5\baselineskip\hcenter{\box\helpbox}{$#2$}\fi%
\def\helpmacro{#3}%
\ifx \helpmacro\empty%
\else\vskip -#5\baselineskip\hcenter{\box\helpbox}{$#3$}\fi%
\def\helpmacro{#4}%
\ifx \helpmacro\empty%
\else\vskip -#5\baselineskip\hcenter{\box\helpbox}{$#4$}\fi%
}}%
\def\relop#1#2#3{\mathrel{\mathop{\kern 0pt #2}\limits^{#1}_{#3}}}%
\def\rahmen#1#2{
 $\vcenter{\hrule%
       \hbox%
        {\vrule%
         \hskip#1%
         \vbox{\vskip#1\relax%
               \hbox{#2}%
               \vskip#1}%
         \hskip#1%
         \vrule}%
       \hrule}$}%
%
%
\def \phantomset#1#2{%
      \setbox\helpbox=\hbox{#1}%
      \hskip 0pt {\hbox to\wd\helpbox{\hss #2\hss}}%
}%
\def\mphantomset #1#2{%
      \setbox\helpbox=\hbox{$#1$}%
      \hskip 0pt {\hbox to\wd\helpbox{\hss $#2$\hss}}%
}%
%
%
%
%
%
%
\def \zitat #1{#1}%
%
%
\def \lit #1#2#3#4{%
\item{[#1]}#2; {\it #3}; {#4}%
}%
%
%
%
\def \terminate{\closeout\reffile \par \vfill \supereject \end}%
\def \terminatebye{\fi\terminate}%
\edef \bye {}%
\def \bye{%
\ifnum \inputdepth >0%
\endinput%
\else%
 \terminatebye%
\fi%
}%
%
%
%
\abovedisplayshortskip=0pt%
\abovedisplayskip=6pt%
\belowdisplayshortskip=0pt%
\belowdisplayskip=6pt%
%
%
\newread \readinfile%
\newread\instream%
\newwrite \reffile%
\newcount\inputdepth%
\inputdepth 0%
\def\readfilepath{}%
\def\readfilesuffix{.tex}%
\def\helpreadfile#1{
  \openin \readinfile=#1\relax%
  \ifeof\readinfile%
    \edef\lokaL{\NOFLAG}%
  \else
    \edef\lokaL{\YESFLAG}%
    \advance\inputdepth 1\relax%
    \input #1\relax%
    \advance\inputdepth -1\relax%
  \fi%
\closein \readinfile%
}%
\def \readfile#1{
  \edef\filenamE{#1}%
  \searchsubstring{\filenamE}{/}%
  \ifnum\RESULT=1%
    \edef\filenamE{\behind}%
    \edef\readfilepath{\infront /}%
  \else\fi%
  \loop%
    \searchsubstring{\filenamE}{/}%
  \ifnum\RESULT=1%
    \edef\filenamE{\behind}%
    \edef\readfilepath{\readfilepath\infront /}%
  \repeat\relax
  \searchsubstring{\filenamE}{.}%
  \ifnum\RESULT=1%
    \edef\filenamE{\infront}%
    \edef\readfilesuffix{.\behind}%
  \else\fi%
  \input \readfilepath\filenamE%
  \ifx\lokaL\NOFLAG%
    \input \filenamE%
  \else\fi%
}%
\def\searchsubstring#1#2{
  \edef\strinG{#1}%
  \edef\substrinG{#2}%
  \ifx\empty\substrinG
    \RESULT 0\relax%
  \else%
    \catcode`\& = 11\relax%
    \edef\dummY{&}
    \edef\strinG{#1\dummY}
    \edef\giveiT{{\strinG}{#2}}%
    \expandafter\helpsearchsubstring\giveiT%
    \ifx\strinG\resultfirstparT
      \RESULT -1\relax%
      \edef\infront{\relax}%
      \edef\behind{\relax}%
    \else
      \RESULT 1\relax%
      \edef\infront{\resultfirstparT}%
      \edef\giveiT{{\resultsecondparT}{\dummY #2}}%
      \expandafter\helpsearchsubstring\giveiT%
      \edef\behind{\resultfirstparT}%
    \fi%
    \catcode`\& = 4\relax%
  \fi%
}%
\def \helpsearchsubstring#1#2{
  \def\TestiT ##1#2##2\Delimiter{%
    \edef\resultfirstparT{##1}%
    \edef\resultsecondparT{##2}%
  }%
  \TestiT#1#2\Delimiter%
}%
\def\replacesubstring #1#2#3{
  \searchsubstring{#1}{#2}%
  \ifnum\RESULT = 1%
    \edef\replaceresult{\infront #3\behind}%
  \else%
    \edef\replaceresult{#1}%
  \fi%
}%
\def\replaceallsubstring #1#2#3{
  \edef\replaceallresult{#1}%
  \edef\lokaL{\NOFLAG}%
  \loop%
    \replacesubstring{\replaceallresult}{#2}{#3}%
    \ifx\replaceallresult\replaceresult%
      \edef\lokaL{\YESFLAG}%
    \else\fi%
    \edef\replaceallresult{\replaceresult}%
  \ifx\lokaL\NOFLAG%
  \repeat\relax%
}%
\def\firstletter#1{
  \edef\lokaLL{{#1}{ }}%
  \expandafter\helpsearchsubstring\lokaLL%
  \if\comparE{ \resultsecondparT}{#1 }%
    \edef\result{ }
  \else
    \edef\lokaLL {#1}%
    \edef\lokaL##1{##1A\lokaLL B}%
    \edef\lokaL{{\expandafter\lokaL \lokaLL}{A\lokaLL B}}%
    \expandafter\helpsearchsubstring\lokaL%
    \let\result\resultfirstparT%
  \fi
}%
\def \foreachletter#1#2{
  \if\comparE{#1}{}%
  \else%
    \firstletter{#1}%
    \searchsubstring{#1}{\result}%
    #2%
    \foreachletter{\behind}{#2}%
  \fi%
}%
\def\stringlength#1{
  \VARX 0\relax
  \foreachletter{#1}{\advance\VARX +1\relax}
}%
%
%
%
%
%
%
%
\def\looP#1\repeaT{\def\bodY{#1}\iteratE}
\def\iteratE{\bodY\let\nexT\iteratE
             \else\let\nexT\relax\fi\nexT}
\let\repeaT=\fi%
\def\loOP#1\repeAT{\def\boDY{#1}\iteraTE}
\def\iteraTE{\boDY\let\neXT\iteraTE
             \else\let\neXT\relax\fi\neXT}%
\let\repeAT=\fi%
%
%
%
%
%
%
%
%
%
%
%
%
\def\loopcommandwithcomma#1\repcommandwithcomma{%
\def\bodycommandwithcomma{#1}\iteratecommandwithcomma}%
\def\iteratecommandwithcomma{%
\bodycommandwithcomma\let\nextcommandwithcomma\iteratecommandwithcomma%
             \else\let\nextcommandwithcomma\relax\fi\nextcommandwithcomma}%
\let\repcommandwithcomma=\fi%
\def \commandwithcomma #1#2{%
  \edef\receiveD{#2}
  \searchsubstring{\receiveD}{,}%
  \ifnum\RESULT<1%
    \edef\argumentcollecT{{\receiveD}}%
  \else%
    \edef\argumentcollecT{}%
    \loopcommandwithcomma%
      \edef\argumentcollecT{\argumentcollecT{\infront}}%
      \edef\argumentremain{\behind}%
      \searchsubstring{\argumentremain}{,}%
    \ifnum \RESULT = 1%
    \repcommandwithcomma%
    \edef\argumentcollecT{\argumentcollecT{\argumentremain}}%
  \fi%
  \expandafter #1\argumentcollecT%
}%
%
%
%
%
%
%
\def\loopmulticommand#1\repmulticommand{%
\def\bodymulticommand{#1}\iteratemulticommand}%
\def\iteratemulticommand{%
\bodymulticommand\let\nextmulticommand\iteratemulticommand%
             \else\let\nextmulticommand\relax\fi\nextmulticommand}%
\let\repmulticommand=\fi%
\def\multicommand #1#2{%
  \edef\lisT{#2}
  \loopmulticommand%
    \searchsubstring{\lisT}{;}%
  \ifnum\RESULT=1%
    \edef\lisT{\behind}%
    \commandwithcomma {#1}{\infront}%
  \repmulticommand\relax%
  \commandwithcomma {#1}{\lisT}%
}%
%
%
%
%
%
%
%
\def\comparedepth #1#2{
  \VARX #1\relax
  \if\comparE{#2}{}
    \VARY 1\relax%
  \else%
    \VARY #2\relax%
  \fi%
  \ifnum\VARX < -\VARX%
   \VARX-\VARX\relax%
  \else\fi
  \edef\result{\NOFLAG}%
  \ifnum\VARX < \VARY%
    #1 \VARY\relax%
    \edef\result{\ONFLAG}%
  \else\fi%
  \ifnum\VARX < -\VARY%
    #1 \VARY\relax
    \edef\result{\OFFFLAG}%
  \else\fi%
}%
%
%
%
%
%
%
%
%
%
%
%
%
\newcount\showlabeldepth%
\showlabeldepth 0\relax%
\edef\showlabelflag{\OFFFLAG}%
\edef \citeext {aux}%
\def\showlabel #1{
  \comparedepth{\showlabeldepth}{#1}%
  \ifx\result\NOFLAG%
  \else%
    \edef\showlabelflag{\result}%
  \fi%
}%
\def\showreferences{%
  \if\comparE{\testshowreferences}{\OFFFLAG}%
  \else%
    \def\testshowreferences{on}%
  \fi%
}%
%
%
%
%
%
\def\compareshowlabel{on}%
\def \auxpath #1{%
\closeout \reffile%
\edef \auxfile {#1\jobname.\citeext}%
\openin \instream=\auxfile\relax%
\ifeof\instream%
\else \input \auxfile\relax%
\fi%
\closein \instream%
\openout \reffile=\auxfile\relax%
}%
\edef \YES {yes}%
\def \auxof #1#2{
  \expandafter\ifx \csname aux#2loaded\endcsname\YES%
  \else%
    \expandafter\edef\csname aux#2loaded\endcsname{yes}\relax%
    \openin \instream=#1#2.\citeext\relax%
    \ifeof\instream%
      \message { FILE   #1#2.\citeext   NOT FOUND}%
      \if\comparE{#1}{}
      \else
        \closein \instream%
        \openin \instream=#2.\citeext\relax%
        \ifeof\instream%
          \message{---------------------------------------------------------------------------}%
          \message { FILE   #2.\citeext   NOT FOUND}%
          \message{---------------------------------------------------------------------------}%
        \else
          \message{// FILE}%
          \input #2.\citeext\relax%
          \message {LOADED //}%
        \fi
      \fi
    \else
      \message{// FILE}%
      \input #1#2.\citeext\relax%
      \message {LOADED //}%
    \fi%
    \closein \instream%
  \fi%
}%
\def \setref #1#2#3{%
  \expandafter \edef\csname R#1R\endcsname {{#2}{#3}{}}%
}%
%
%
%
%
\catcode`\_ = 14%
\catcode`\% = 12_
\edef\percentsigN{
\catcode`\% = 14_
\catcode`\_ = 8%
\def\loadmainauxfile{
  \edef\auxfile{\jobname.\citeext}%
  \openin \instream=\auxfile\relax%
  \ifeof\instream%
  \else%
    \input \auxfile\relax%
  \fi%
  \closein \instream%
  \openout \reffile =\auxfile\relax%
  \def \setref ##1##2##3{%
    \edef \refmacro{\csname R##1R\endcsname}
    \expandafter\ifx \refmacro\relax%
      \expandafter \edef\csname R##1R\endcsname {{##2}{##3}{\auxsource}}%
    \else%
    \fi%
  }%
  \edef \writelabel{%
    \write\reffile {\edef\string\auxsource {\jobname}\percentsigN}%
  }%
  \writelabel
}
\def\label #1#2{
  \edef\parameterone{#1}%
  \ifx \parameterone\empty
  \else%
    \ifx \showlabelflag\compareshowlabel%
      \rahmen{1pt}{{\tt \small{\parameterone}}}%
    \else%
    \fi%
    \expandafter\edef\csname R#1R\endcsname{{#2}{\folio}{}}
    \edef \writelabel{%
    \write\reffile {\string\setref {#1}{#2}{\noexpand\folio}\percentsigN}}%
    \writelabel
  \fi%
}%

%
\def \writereference#1#2#3#4{
  \ifx #3\empty%
  \else%
    #3:%
  \fi%
  \ifx #4\empty%
    #1%
  \else%
    #2%
  \fi%
}%
\def \ref#1{
  \edef \refmacro{\csname R#1R\endcsname}
  \expandafter\ifx \refmacro\relax
 {\bf #1}%
    \message{---------------------------------------------------------------------------}%
    \message {#1--  UNDEFINIERTE REFERENZ}%
    \message{---------------------------------------------------------------------------}%
    \message {}%
  \else%
    \expandafter\writereference\refmacro{\empty }%
  \fi%
  \ifx \testshowreferences\compareshowlabel%
    \rahmen{1pt}{{\tt \small{#1}}}%
  \else\fi%
}%
\def \pageref#1{
  \edef \refmacro{\csname R#1R\endcsname}%
  \expandafter\ifx \refmacro\relax
    {\bf #1}%
  \else%
    \expandafter\writereference\refmacro{\bigskip }
  \fi%
}%
%
%
\def \gluefiles{\noaux}%
\def\makegoodjobname{
  \edef \lokaL ##1##2##3##4##5##6##7##8{##1##2##3##4##5##6}%
  \edef\goodjobname{\expandafter\lokaL\jobname ........}%
  \searchsubstring{\goodjobname}{.}%
  \edef\goodjobname{\infront}%
}%
%
%
%
%
%
%
%
%
\ifx\diagramloaded\relax\endinput\else\let\diagramloaded\relax\fi%
\chardef\body\catcode`@\catcode`@=11%
\let\q@wlog\wlog\def\wlog#1{} \def\q@del#1\@{}%
\newcount\dpi\dpi300 \newcount\q@fn\q@fn\z@\newcount\q@cn\q@cn255%
\newdimen\q@fct\q@fct10pt \newdimen\q@scr%
\newwrite\q@file\newif\ifq@open\newif\ifrztex%
\def\q@wrtc{\afterassignment\q@wr@c\count@}%
\def\q@wrtd{\afterassignment\q@wr@d\dimen@}%
\def\q@wr@c{\immediate\write\q@file{\the\count@}}%
\def\q@wr@d{\expandafter\q@wr@e\the\dimen@}%
\def\q@cut{\expandafter\q@cu@\the}%
\let\q@wr@a\immediate\let\q@wr@b\write%
{\catcode`p=12\catcode`t=12 \gdef\q@cu@#1pt{#1}%
  \gdef\q@wr@e#1pt{\q@wr@a\q@wr@b\q@file{#1}}}%
\def\q@setl#1#2{\edef\q@tx{\the#1\noexpand\or#2}#1\expandafter{\q@tx}}%
\def\q@div#1#2{\q@scr#2\divide\q@scr2048\divide#1\q@scr\multiply#1 32}%
\newbox\q@diabx\newif\ifq@dia\newcount\q@m\newdimen\q@x\newdimen\q@y%
\newdimen\q@medmuskip%
\def\diagram{\relax\ifq@dia\q@erra\fi\ifq@open\else\global\q@opentrue%
    \immediate\openout\q@file diagram.stf%
    \immediate\write\q@file{\goodjobname}\q@wrtc\dpi\q@wrtc\mag\fi%
  \q@xmin\maxdimen\q@xmax-\maxdimen\q@ymin\maxdimen\q@ymax-\maxdimen%
  \setbox\q@diabx\hbox\bgroup\m@th%
  \setbox\z@\hbox{$\mkern\medmuskip$}\q@medmuskip\wd\z@%
  $\q@diatrue\q@init%
  \ifrztex\q@fct\Gradof\@ktGrad\fi\q@wrtd100\q@fct%
  \let\object\q@@object\let\q@ext\q@@ext\let\arrow\q@@arrow}%
\def\enddiagram{\q@makelist\q@m\q@cnt%
  \loop\ifnum\q@m>\z@\q@list\q@getxy\advance\q@x\q@c\advance\q@y-.25\q@fct%
    \setbox\q@obj\hbox{\unhbox\q@obj\global\setbox\z@\lastbox}\q@put\z@%
    \advance\q@m\m@ne\repeat$\egroup\ifdim\q@xmin=\maxdimen\else%
  \setbox\q@diabx\hbox{\kern-\q@xmin\box\q@diabx\kern\q@xmax}%
  \ht\q@diabx\q@ymax\dp\q@diabx-\q@ymin\ifmmode\vcenter{\box\q@diabx}%
  \else$\vcenter{\box\q@diabx}$\fi\fi\q@wrtc\m@ne}%
\def\q@erra{\errmessage{Don't construct diagrams within diagrams}}%
\def\q@errb{\errmessage{Use \string\freeunits\space before \string\object}}%
\def\q@errc{\errmessage{Don't use \string\gapsadjust\space or \string%
  \small\space in \string\freeunits-mode}}%
\def\q@errd{\errmessage{Indexes may not exceed 100 when%
  \string\gapsadjust\space or \string\small\space is used}}%
\def\q@erre{\errmessage{Target-point is too close to starting-point%
  (minimum is 2pt)}}%
\def\q@errf{\errmessage{Use \string\hgaps, \string\vgaps\space and%
  \string\gapsadjust\space before \string\arrow}}%
\def\q@errg{\errmessage{In this context use positive numbers only}}%
\def\q@errh{\errmessage{Don't use \string\object, \string\arrow\space%
  or \string\gapsadjust\space outside the \string\diagram-environment}}%
\newdimen\hunit\newdimen\vunit\newcount\q@hmax\newcount\q@vmax%
\def\unit{\relax\afterassignment\q@unit\hunit} \def\q@unit{\vunit\hunit}%
\def\hgaps#1{\relax\q@bg@\hunit#1;;\@\q@bgd\q@hgaps\hunit\q@hmax}%
\def\vgaps#1{\relax\q@bg@\vunit#1;;\@\q@bgd\q@vgaps\vunit\q@vmax}%
\def\q@bg@#1{\dimen@ii#1\dimen@\z@\count@\m@ne\toks@{}\q@bga}%
\def\q@bga{\futurelet\next\q@bgb}%
\def\q@bgb{\ifx\next;\let\next\q@del\else\let\next\q@bgc\fi\next}%
\def\q@bgc#1;{\q@setl\toks@{\dimen@\the\dimen@}\afterassignment\q@del%
  \advance\dimen@#1\dimen@ii\@\advance\count@\m@ne\q@bga}%
\def\q@bgd#1#2#3{#3-\count@\edef#1{\relax\noexpand\ifcase\count@%
  \the\toks@\noexpand\else\advance\count@\the\count@\dimen@\the\dimen@%
  \advance\dimen@\count@#2\noexpand\fi}}%
\newdimen\q@a\newdimen\q@b\newdimen\q@c\newdimen\q@d\newdimen\q@e%
\newcount\q@i\newcount\q@j\newbox\q@obj\setbox\q@obj\null%
\newtoks\q@toks\q@toks{} \newtoks\q@@toks\q@@toks{}%
\newcount\q@homax\newcount\q@vomax\q@homax\z@\q@vomax\z@%
\newif\ifq@free\newif\ifq@list\newif\ifq@xlist%
\newtoks\q@ta\newtoks\q@tb\newtoks\q@tc\newcount\q@cnt\q@cnt\z@%
\def\freeunits{\q@freetrue\let\q@ext\q@errc}%
\newdimen\q@tabkern\newdimen\q@tbkern%
\newtoks\q@tcapprox\newtoks\q@taapprox%
\def\q@tabtckind#1#2{{\thinmuskip0mu\thickmuskip0mu\medmuskip18mu%
 \setbox\z@\hbox{$#1$}\dimen@\wd\z@\advance\dimen@.1\p@%
 \setbox\z@\hbox{$#1{}$}%
 \ifdim\wd\z@>\dimen@\setbox\z@\hbox{$#2$}\advance\dimen@\wd\z@%
  \setbox\z@\hbox{$#1#2$}%
  \ifdim\wd\z@>\dimen@\global\q@tabkern\q@medmuskip%
   \global\q@tbkern.5\q@medmuskip\global\q@tcapprox{{}}\fi\fi}}%
\def\q@tatbkind#1#2{{\thinmuskip0mu\thickmuskip0mu\medmuskip18mu%
 \setbox\z@\hbox{$#1$}\dimen@\wd\z@\advance\dimen@.1\p@%
 \setbox\z@\hbox{$#1{}$}%
 \ifdim\dimen@>\wd\z@\setbox\z@\hbox{$#2$}\advance\dimen@\wd\z@%
  \setbox\z@\hbox{$#1#2$}%
  \ifdim\wd\z@>\dimen@\global\advance\q@tbkern.5\q@medmuskip%
   \global\q@taapprox{{}}\fi\fi}}%
\def\q@tbcal{\dimen@.5\wd\z@\advance\dimen@-\q@tbkern%
 \q@e\dimen@\advance\q@e-.8\q@fct\ifdim\q@e<\z@\q@e\z@\fi}%
\def\q@tabcal{\q@c-\wd\z@\advance\q@c\q@tabkern\advance\q@c\dimen@}%
\def\q@@object#1(#2;#3)#4{\relax\let\freeunits\q@errb\advance\q@cnt\@ne%
  \ifq@free\afterassignment\q@del\q@x#2\hunit\@\afterassignment\q@del%
    \q@y#3\vunit\@\q@setl\q@@toks{\q@x\the\q@x\q@y\the\q@y}%
  \else\ifnum#2<\@ne\q@errg\fi\ifnum#3<\@ne\q@errg\fi%
    \q@setl\q@@toks{\q@i#2 \q@j#3 }%
    \ifnum#2>\q@vomax\q@vomax#2\fi\ifnum#3>\q@homax\q@homax#3\fi\fi%
  \q@ta{}\q@tb{}\q@tc{}\setbox\z@\vbox{\ialign{\global\q@tb{##}%
  &\global\q@tc{##}&\global\q@ta\q@tb\global\q@tb\q@tc\global\q@tc{##}%
  &&\global\q@tc\expandafter{\the\q@tc##}\crcr#4\crcr}}%
  \q@tabkern\z@\q@tbkern\z@\q@tcapprox{}\q@taapprox{}%
  \edef\q@test{\the\q@ta}%
  \ifx\q@test\empty\edef\q@test{\the\q@tc}%
   \ifx\q@test\empty%
    \setbox\z@\hbox{$\the\q@tb$}\q@tbcal\q@tabcal%
   \else%
    \q@tabtckind{\the\q@tb}{\the\q@tc}%
    \setbox\z@\hbox{$\the\q@tb\the\q@tcapprox$}\q@tbcal\q@tabcal%
    \setbox\z@\hbox{$\the\q@tb\the\q@tc$}\fi%
  \else%
    \q@tabtckind{\the\q@ta\the\q@tb}{\the\q@tc}%
    \q@tatbkind{\the\q@ta}{\the\q@tb\the\q@tcapprox}%
    \setbox\z@\hbox{$\the\q@taapprox\the\q@tb\the\q@tcapprox$}\q@tbcal%
    \setbox\z@\hbox{$\the\q@ta\the\q@tb\the\q@tcapprox$}\q@tabcal%
    \setbox\z@\hbox{$\the\q@ta\the\q@tb\the\q@tc$}\fi%
    \q@a\wd\z@\advance\q@a\q@c\q@b\ht\z@\advance\q@b-.25\q@fct%
    \q@d-\dp\z@\advance\q@d-.25\q@fct \q@setl\q@toks{\q@a\the\q@a%
  \q@b\the\q@b\q@c\the\q@c\q@d\the\q@d\q@e\the\q@e}%
  \q@listfalse\q@xlistfalse\setbox\q@obj\hbox{\unhbox\q@obj\box\z@}}%
\def\gapsadjust{\hgapsadjust\vgapsadjust}%
\def\hgapsadjust{\q@adj\q@hmax\q@homax\q@hgaps\hunit-\q@j\q@a>\q@c}%
\def\vgapsadjust{\q@adj\q@vmax\q@vomax\q@vgaps\vunit{}\q@i\q@d<\q@b}%
\def\q@@ext#1#2#3#4#5{\relax\ifnum#1>100\q@errd\fi\q@makelist%
  \bgroup\dimendef\dimen@0\dimendef\dimen@i250\dimendef\q@e251%
  \dimendef\q@a252\dimendef\q@b253\dimendef\q@c254\dimendef\q@d255%
  \count@\z@\loop\ifnum\count@<#1\advance\count@101%
    \dimen\count@\z@\advance\count@-100\dimen\count@\z@\repeat%
  \q@m\z@\loop\ifnum\q@m<\q@cnt\advance\q@m\@ne\q@list\q@@list%
    \ifdim#3#4\dimen#2\dimen#2#3\fi\advance#2100%
    \ifdim\dimen#2#4#5\dimen#2#5\fi\repeat\toks@{}}%
\def\q@adj#1#2#3#4#5#6#7#8#9{\q@ext#2#6#7#8#9\ifnum#1<#2#1#2\fi%
  \dimen@\z@\dimen@i\z@\q@m\@ne \loop\ifnum\q@m<#1\q@setl\toks@{\dimen@%
    \the\dimen@}\ifnum\q@m<#2\advance\dimen@i-#5\dimen\q@m\advance\q@m101%
    \advance\dimen@i#5\dimen\q@m\advance\q@m-100 \else\advance\q@m\@ne \fi%
    \count@\q@m#3\advance\dimen@\dimen@i\repeat\global\toks@\toks@%
  \global\dimen@\dimen@\global\count@-#1\egroup\q@bgd#3#4#1}%
\def\q@crt#1#2#3#4#5#6{\q@ext#5#2#3#6#4\q@m\z@\loop\ifnum\q@m<\q@cnt%
  \advance\q@m\@ne \q@@list\dimen@\dimen#2\advance#2100%
    \q@setl\toks@{#3\the\dimen@#4\the\dimen#2}\repeat\global\toks@\toks@%
  \egroup\edef#1{\noexpand\ifcase\q@m\the\toks@\noexpand\fi}}%
\def\q@makelist{\ifq@list\else\q@listtrue%
  \edef\q@list{\noexpand\ifcase\q@m\the\q@toks\noexpand\fi}%
  \edef\q@@list{\noexpand\ifcase\q@m\the\q@@toks\noexpand\fi}\fi}%
\def\q@getxy{\q@@list\ifq@free\else%
  \count@\q@j\q@hgaps\q@x\dimen@\count@\q@i\q@vgaps\q@y-\dimen@\fi}%
\newif\ifq@frcor\let\fromcorner\q@frcortrue%
\newif\ifq@tocor\let\tocorner\q@tocortrue%
\newif\ifq@swap\newif\ifq@arc\newif\ifq@swname%
\newif\ifq@sml\let\small\q@smltrue%
\newbox\q@nbox\newcount\q@ncnt\newtoks\q@nlist%
\newcount\q@equl\newcount\q@brk\newcount\q@head\newcount\q@epi%
\newcount\q@bth\newcount\q@map\newcount\q@mono\newcount\q@iso%
\newcount\q@incl\newcount\q@ucnt\newtoks\q@ulist%
\newdimen\q@xa\newdimen\q@ya\newdimen\q@xb\newdimen\q@yb%
\newdimen\q@dxa\newdimen\q@dya\newdimen\q@dxb\newdimen\q@dyb%
\newdimen\q@k\newdimen\q@l\newcount\q@ma\newcount\q@mb%
\newdimen\q@dx\newdimen\q@dy \newdimen\q@h%
\newdimen\nameskip\nameskip2pt%
\newdimen\brokenline\brokenline3pt%
\newdimen\brokengap\brokengap2pt%
\newdimen\displace\newdimen\fromskip\newdimen\toskip%
\def\fnull{\fromskip\z@}\def\tnull{\toskip\z@}%
\newdimen\q@skip\newdimen\q@arrsk%
\def\arrowskip{\relax\afterassignment\q@sk@p\q@arrsk}%
\def\q@sk@p{\toskip\q@arrsk\fromskip\q@arrsk} \arrowskip5pt%
\def\q@init{\q@frcorfalse\q@tocorfalse\q@swapfalse\q@arcfalse%
  \q@smlfalse\q@swnamefalse\q@equl\z@\q@brk\z@\q@head8%
  \q@epi\z@\q@bth\z@\q@map\z@\q@incl\z@\q@iso\z@\q@mono\z@%
  \q@ucnt\z@\q@ulist{}\q@ncnt\z@\q@nlist{}\q@sk@p\displace\z@}%
\def\arc{\relax\q@arctrue\q@h} \def\equal{\relax\q@equl2 }%
\def\broken{\relax\q@brk4 } \def\nohead{\relax\q@head\z@}%
\def\epi{\relax\q@epi16 } \def\both{\relax\q@bth32 }%
\def\map{\relax\q@map64 } \let\swap\q@swaptrue%
\def\incl{\relax\ifq@swap\q@incl256 \else\q@incl128 \fi\q@swapfalse}%
\def\iso{\relax\ifq@swap\q@iso1024 \else\q@iso512 \fi\q@swapfalse}%
\def\mono{\relax\q@mono2048 }%
\def\under{\relax\afterassignment\q@und\count@}%
\def\q@und{\advance\q@ucnt\@ne \edef\q@tx{\the\q@ulist\noexpand%
  \q@wrtc\the\count@}\q@ulist\expandafter{\q@tx}}%
\def\name{\relax\futurelet\next\q@name}%
\def\q@name{\ifx\next(\let\next\q@n@me%
  \else\let\next\relax\def\next{\q@n@me(\z@;\z@)}\fi\next}%
\def\q@n@me(#1;#2)#3{\advance\q@ncnt\@ne%
  \q@setl\q@nlist{\q@dx#1\q@dy#2\ifq@swap\q@d-\q@d\fi}\q@swapfalse%
  \setbox\q@nbox\hbox{\unhbox\q@nbox\hbox{$\scriptstyle#3$}}}%
\def\q@@arrow{\relax\let\q@bgd\q@errf\futurelet\next\q@arrow}%
\def\q@arrow{\ifx\next(\let\next\q@@rrow\else%
  \let\next\relax\def\next{\q@@rrow(\z@;\z@)(\z@;\z@)}\fi\next}%
\def\q@@rrow(#1;#2)#3(#4;#5){\q@dxa#1\q@dya#2\q@dxb#4\q@dyb#5%
  \afterassignment\q@@@r@w\q@ma}%
\def\q@@@r@w{\afterassignment\q@@rr@@\q@mb}%
\def\q@@rr@@{\ifnum\q@ma>\z@\ifnum\q@mb>\z@%
  \ifnum\q@ma>\q@cnt\else\ifnum\q@mb>\q@cnt\else\q@makelist%
  \q@m\q@ma\q@getxy\q@list\q@xa\q@x\q@ya\q@y \dimen@i\q@e%
  \q@m\q@mb\q@getxy\q@list\q@xb\q@x\q@yb\q@y%
  \ifdim\q@xa<\q@xb\q@e-\q@e\else\ifdim\q@xa>\q@xb\dimen@i-\dimen@i%
    \else\q@e\z@\dimen@i\z@\fi\fi%
  \ifq@frcor\advance\q@xa\dimen@i\fi\ifq@tocor\advance\q@xb\q@e\fi%
  \advance\q@xa\q@dxa\advance\q@ya\q@dya%
  \advance\q@xb\q@dxb\advance\q@yb\q@dyb%
  \ifdim\q@xa>\q@xb\q@swnametrue%
    \ifnum\q@incl>\z@\count@-\q@incl\q@incl384 \advance\q@incl\count@\fi%
    \ifnum\q@iso>\z@\count@-\q@iso\q@iso1536 \advance\q@iso\count@\fi\fi%
  \ifq@arc\count@\@ne \else\count@\z@\fi%
  \advance\count@\q@equl\advance\count@\q@brk\advance\count@\q@head%
  \advance\count@\q@epi\advance\count@\q@bth\advance\count@\q@map%
  \advance\count@\q@incl\advance\count@\q@iso\advance\count@\q@mono%
  \q@wrtc\count@\q@wrtc\q@ucnt\the\q@ulist%
  \ifnum\q@brk>\z@\q@wrtd\brokenline\q@wrtd\brokengap\fi%
  \q@x\q@xb\advance\q@x-\q@xa\q@y\q@yb\advance\q@y-\q@ya%
  \q@k\q@x\q@l\q@y\divide\q@x4\divide\q@y4%
  \dimen@ii\q@cut\q@x\q@x\advance\dimen@ii\q@cut\q@y\q@y%
  \ifdim\dimen@ii<.2pt\q@erre\fi%
  \dimen@\dimen@ii\advance\dimen@-1600pt%
  \divide\dimen@80\advance\dimen@40pt%
  \dimen@i\dimen@\dimen@\dimen@ii\q@div\dimen@\dimen@i%
  \advance\dimen@\dimen@i\divide\dimen@2%
  \dimen@i\dimen@\dimen@\dimen@ii\q@div\dimen@\dimen@i%
  \advance\dimen@\dimen@i\divide\dimen@2%
  \dimen@i\dimen@\dimen@\dimen@ii\q@div\dimen@\dimen@i%
  \advance\dimen@\dimen@i\divide\dimen@2%
  \dimen@i\dimen@\dimen@\dimen@ii\q@div\dimen@\dimen@i%
  \advance\dimen@\dimen@i\divide\dimen@2%
  \dimen@i\dimen@\dimen@\dimen@ii\q@div\dimen@\dimen@i%
  \advance\dimen@\dimen@i\divide\dimen@2%
  \dimen@i\dimen@\dimen@\dimen@ii\q@div\dimen@\dimen@i%
  \advance\dimen@\dimen@i\divide\dimen@2%
  \dimen@i\dimen@\dimen@\dimen@ii\q@div\dimen@\dimen@i%
  \advance\dimen@\dimen@i\divide\dimen@2%
  \dimen@i\dimen@\dimen@\dimen@ii\q@div\dimen@\dimen@i%
  \advance\dimen@\dimen@i\divide\dimen@2%
  \dimen@i\dimen@\dimen@\dimen@ii\q@div\dimen@\dimen@i%
  \advance\dimen@\dimen@i%
  \dimen@2\dimen@\q@div\q@k\dimen@\q@div\q@l\dimen@%
  \q@skip\fromskip\ifnum\q@mono>\z@\advance\q@skip.204\q@fct\else%
  \ifnum\q@incl>\z@\advance\q@skip.13\q@fct\fi\fi%
  \q@corr\q@xa\q@ya\q@ma\q@skip\q@k-\q@k\q@l-\q@l%
  \q@corr\q@xb\q@yb\q@mb\toskip%
  \dimen@\q@k\q@k\q@l\q@l-\dimen@\displace-\displace%
  \advance\q@xa\q@cut\displace\q@k\advance\q@ya\q@cut\displace\q@l%
  \advance\q@xb\q@cut\displace\q@k\advance\q@yb\q@cut\displace\q@l%
  \q@wrtd\q@xa\q@wrtd\q@ya\q@wrtd\q@xb\q@wrtd\q@yb%
  \ifq@arc\q@wrtd\q@h\q@wrtd\q@skip\q@wrtd\toskip\q@h-\q@h\fi%
  \global\advance\q@cn\@ne \ifnum\q@cn=256%
    \global\advance\q@fn\@ne%
    \global\font\q@font=\jobname\ifcase\q@fn\or\or _b\or _c\or _d%
      \or _e\or _f\or _g\or _h\or _i\or _j\or _k\or _l\or _m\or _n\or _o\or _p%
      \or _q\or _r\or _s\or _t\fi%
    \global\q@cn\z@\fi\hbox{\q@font\char\q@cn}%
  \advance\q@xa\q@xb\divide\q@xa2\advance\q@ya\q@yb\divide\q@ya2%
  \ifq@arc\advance\q@xa\q@cut\q@h\q@k\advance\q@ya\q@cut\q@h\q@l\fi%
  \edef\q@nlst{\noexpand\ifcase\q@ncnt\the\q@nlist\noexpand\fi}%
  \loop\ifnum\q@ncnt>\z@%
    \setbox\q@nbox\hbox{\unhbox\q@nbox\global\setbox\z@\lastbox}%
    \q@a\wd\z@\divide\q@a2\q@b\ht\z@\advance\q@b\dp\z@\divide\q@b2%
    \q@c-\q@b\advance\q@c\dp\z@%
    \ifdim\q@k>\z@\ifdim\q@l<\z@\q@b-\q@b\fi%
      \else\ifdim\q@l>\z@\q@b-\q@b\fi\fi%
    \q@d\q@cut\q@a\q@k\advance\q@d\q@cut\q@b\q@l%
    \ifdim\q@d<\z@\q@d-\q@d\fi\advance\q@d\nameskip%
    \ifq@swname\q@d-\q@d\fi\q@nlst\advance\q@d\q@dy \q@x\q@xa\q@y\q@ya%
    \advance\q@x\q@cut\q@d\q@k\advance\q@y\q@cut\q@d\q@l%
    \advance\q@x\q@cut\q@dx\q@l\advance\q@y-\q@cut\q@dx\q@k%
    \advance\q@x-\q@a\advance\q@y\q@c\q@put\z@\advance\q@ncnt\m@ne\repeat%
  \q@wrtc1711 \q@init\fi\fi\fi\fi}%
\def\q@corr#1#2#3#4{\q@m#3\q@getxy\ifq@sml\ifq@xlist\else%
    \q@xlisttrue\q@crt\q@aclist\q@j\q@a\q@c\q@homax>\q@crt\q@bdlist\q@i%
      \q@d\q@b\q@vomax<\fi\q@aclist\q@bdlist\else\q@list\fi%
  \dimen@#1\advance\dimen@-\q@x\advance\q@a-\dimen@\advance\q@c-\dimen@%
  \dimen@#2\advance\dimen@-\q@y\advance\q@b-\dimen@\advance\q@d-\dimen@%
  \ifq@arc\q@wrtd\q@a\q@wrtd\q@b\q@wrtd\q@c\q@wrtd\q@d\else%
  \ifdim\q@a>\z@\ifdim\q@b>\z@\ifdim\q@c<\z@\ifdim\q@d<\z@%
    \ifdim\q@l>\z@\ifdim\q@k>\z@\q@c@rr\q@a\q@b>\else\q@c@rr\q@c\q@b<\fi%
    \else\ifdim\q@k>\z@\q@c@rr\q@a\q@d<\else\q@c@rr\q@c\q@d>\fi\fi%
  \advance#1\dimen@\advance#2\dimen@i\fi\fi\fi\fi%
  \advance#1\q@cut#4\q@k\advance#2\q@cut#4\q@l\fi}%
\def\q@c@rr#1#2#3{\dimen@i\q@cut#1\q@l\dimen@\q@cut#2\q@k%
  \ifdim\dimen@i#3\dimen@\q@div\dimen@\q@l\dimen@i#2%
    \else\q@div\dimen@i\q@k\dimen@#1\fi}%
\newdimen\q@xmax\newdimen\q@xmin\newdimen\q@ymax\newdimen\q@ymin%
\def\q@put#1{\ifdim\q@x<\q@xmin\global\q@xmin\q@x\fi%
  \dimen@\q@x\advance\dimen@\wd#1%
    \ifdim\dimen@>\q@xmax\global\q@xmax\dimen@\fi%
  \dimen@\q@y\advance\dimen@-\dp#1%
    \ifdim\dimen@<\q@ymin\global\q@ymin\dimen@\fi%
  \dimen@\q@y\advance\dimen@\ht#1%
    \ifdim\dimen@>\q@ymax\global\q@ymax\dimen@\fi%
  \rlap{\kern\q@x\raise\q@y\box#1}}%
\let\object\q@errh\let\q@ext\q@errh\let\arrow\q@errh%
\hgaps{}\vgaps{}\unit1cm \let\wlog\q@wlog\catcode`@\body 
\makegoodjobname%
\loadmainauxfile%
\pageno 1\relax%
\fi%

\def \writereference#1#2#3#4{
  \ifx #3\empty%
  \else%
    \zitat{[Tr], }
  \fi%
  \ifx #4\empty%
    #1%
  \else%
    #2%
  \fi%
}%
\edef \auxsource{CUTSA}%
\setref{def gheir}{(2.1)}{2}%
\setref{explain heir}{(2.2)}{2}%
\setref{general test for heir}{(2.3)}{2}%
\setref{f(p)}{(2.4)}{2}%
\setref{heirs of 1-types}{(2.5)}{3}%
\setref{characterization of dense types}{(3.1)}{4}%
\setref{dense type}{(3.2)}{4}%
\setref{qe dense}{(3.3)}{4}%
\setref{invariance group}{(3.4)}{5}%
\setref{multiplicative invariance group}{(3.5)}{5}%
\setref{connect hatp and tildep}{(3.6)}{5}%
\setref{signature makes sense}{(3.7)}{6}%
\setref{signature}{(3.8)}{6}%
\setref{prop of sign}{(3.9)}{6}%
\setref{sim}{(3.10)}{7}%
\setref{general reduction}{(3.11)}{7}%
\setref{reformulation of alternative}{(3.12)}{7}%
\setref{signature alternative}{(3.13)}{8}%
\setref{test if only one side}{(3.14)}{8}%
\setref{helppphat}{(3.15)}{8}%
\setref{pphat}{(3.16)}{8}%
\setref{T-convex}{(4.1)}{9}%
\setref{T-convex condition}{(4.2)}{9}%
\setref{describe T-convex}{(4.3)}{9}%
\setref{structure on residue field}{(4.4)}{9}%
\setref{model theory of $T$-convex}{(4.5)}{9}%
\setref{explicit description}{(4.6)}{10}%
\setref{quote extended Marker-Steinhorn}{(4.7)}{10}%
\setref{invariance valuation ring}{(5.1)}{10}%
\setref{charact of definable in Tconvex}{(5.2)}{10}%
\setref{signature alternative for groups}{(5.4)}{12}%
\setref{heirs for groups}{(5.5)}{12}%
\setref{mc groups}{(5.6)}{12}%
\setref{drain lemma}{(6.1)}{12}%
\setref{obtain sign 0}{(6.2)}{13}%
\setref{extend to sign 0}{(6.3)}{13}%
\setref{chains}{(6.4)}{14}%
\setref{obtain both signature 0}{(6.5)}{14}%
\setref{obtain arbitrary signs}{(6.6)}{15}%
\setref{locate derivative}{(7.1)}{16}%
\setref{signature alternatives}{(7.2)}{16}%
\setref{pphatVp}{(7.3)}{17}%
\setref{mc}{(7.4)}{18}%
\headlineshape {}{off}{double}
\title {}{}{\zl{The elementary theory of Dedekind cuts in polynomially bounded structures}}{}

{\zl{\it Marcus Tressl.}} 
{\zl{{\it NWF-I Mathematik, 93040 Universit\"at Regensburg, Germany}}}
\zl{\it e-mail: marcus.tressl@mathematik.uni-regensburg.de}
\footnote{}{2000 Mathematics Subject Classification: Primary 03C64;
Secondary 03C10}
\showlabel {0}
\longversion{0}{1}
\medskip
\medskip
\abstract{
Let $M$ be a polynomially bounded, o-minimal structure with archi\-medean prime model, for example if $M$ is a real closed field.
Let $C$ be a convex and unbounded subset of $M$.
We determine the first order theory of the structure $M$ expanded
by the set $C$.
We do this also over any given set of parameters from $M$,
which yields a description of all subsets of $M^n$, definable
in the expanded structure.}

\title{}{*}{Introduction}{}

This paper is a sequel to \zitat{[Tr]}, where
we began the model theoretic  study of Dedekind cuts
of o-minimal expansions of fields. Before explaining what we do here,
we recall some terminology from \zitat{[Tr]}.\einzug
If $X$ is a totally ordered set, then a (Dedekind) cut $p$ of $X$ is a tuple
$p=(p^L,p^R)$ with $X=p^L\cup p^R$ and $p^L<p^R$. If $Y\subseteq
X$ then $Y^+$ denotes the cut $p$ of $X$ with $p^R=\{ x\in X\st
x>Y\} $. $Y^+$ is called the upper edge of $Y$. Similarly the
lower edge $Y^-$ of $Y$ is defined.\medskip

We fix an o-minimal expansion $T$ of the theory of real closed fields
in a language $\SL $. If $M$ is a model of $T$ and $p$ is a cut of
(the underlying set of) $M$, then by the model theoretic
properties of $p$ we understand the model theoretic
properties of the structure $M$ expanded by the set $p^L$
in the language $\SL (\SD )$ extending $\SL $, which has an additional
unary predicate $\SD $ - interpreted as $p^L$.
We write $(M,p^L)$ for this expansion of $M$.
\einzug
Our aim here is to determine the full theory of the structure $(M,p^L)$
in the language $\SL (\SD )$ relative $T$ and to give a description
of the subsets of $M^n$, definable in $(M,p^L)$ relative $M$.
By ``relative $T$" and ``relative $M$" we mean that the theory $T$
and the subsets of $M^n$ definable in $M$ are assumed to be known.
\einzug
By basic model theory, this problem amounts
to find the theory of the structure $(M,p^L)$
in the language $\SL (\SD )$ over a given set $A$ of parameters.
\medskip
We can do this for all cuts in the case where $T$ is polynomially bounded
with archimedean prime model
(c.f \ref {def poly bound} below), in particular in the case
of pure real closed fields.
The main result is Theorem \ref {general main} below, which is of technical
nature. For the moment, we describe what we get from this result by saying what the
subsets of $M^n$, definable in $(M,p^L)$, are.
In order to speak about these sets, we first have to recall some
invariants of a cut $p$ from \zitat {[Tr]}. The o-minimal
terminology is taken from \zitat{[vdD1]}.
\medskip
\medskip

\Definition{}{}
Let $p$ be a cut of an ordered abelian group $K$.
The convex subgroup
$$G(p):=\{ a\in K\st a+p=p\} $$
of $K$ is called the \notion {invariance group}
of $p$ (here $a+p:=(a+p^L,a+p^R)$).
The cut $G(p)^+$ is denoted by $\hat p$.
\smallskip
\einzug
Now let $K$ be an ordered field. Then $G(p)$
denotes the invariance group of $p$
with respect to $(K,+,\leq )$ and $G^*(p)$ denotes
the invariance group of $\a p\a $ with respect to
$(K^{>0},\cdot ,\leq )$, hence
$G^*(p)=\{ a\in K\st a\mal p=p\} $.\einzug
Moreover, the convex valuation ring
$$V(p):=\{ a\in K\st a\mal G(p)\subseteq G(p)\} $$
of $K$ is called the \notion {invariance ring} of $p$.
Note that the group of positive units of $V(p)$ is the multiplicative invariance group of $\hat p$.\einzug
If $X$ is a symmetric subset of $K$ then we write
$G(X)$ and $V(X)$ for $G(X^+)$ and $V(X^+)$, respectively.
If $s\not \in K$ is from an ordered field extension of $K$
then we write $G(s/K)$, $G^*(s/K)$ and $V(s/K)$ for the invariance groups
and the invariance ring of the cut induced by $s$ on $K$.
\medskip
By \ref {connect hatp and tildep}
if $p>\hat p$, then there is some
$c\in K$ such that
$$G^*(p)=c\mal G(p)+1\ (=\{ c\mal a+1\st a\in G(p)\} )$$
\medskip
\Definition {signature}
Let $K$ be a divisible ordered abelian group and let $p$ be a cut of
$K$. We define the \notion {signature} of $p$ as
$$\sign p:=\cases{1 & if there is a convex subgroup $G$ of $K$ and some
$a\in K$ with $p=a+G^+$\cr
-1 & if there is a convex subgroup $G$ of $K$ and some $a\in K$ with
$p=a-G^+$\cr
0 & otherwise}$$

Since $K$ is divisible we can not have $a+G^+=b-H^+$
for $a,b\in K$ and convex subgroups $G,H$ of $K$.
Hence the signature is well defined.\einzug
If $K$ is a real closed field, then $\sign ^*p$
denotes the signature of $\a p\a $ with respect
to $(K,\cdot ,\leq )$.
\medskip

\Definition{}{}
Let $p$ be a cut of a divisible ordered abelian group
$M$. If $\sign p\neq 0$, then $p$ is an edge of the nonempty set
$Z(p):=\{ a\in M\st a+\hat p=p\text { or }a-\hat p=p\} $
and we denote the other edge by $z(p)$.
If $\sign p=0$ we define $z(p):=p$.\einzug
If $G$ is a convex subgroup of a real closed field $M$,
we write $Z^*(G)$ and $Z^*(G^+)$ for the
set $\{ a\in M\st a\mal V(G)^+=G^+\text { or }a\mal \Dm (V(G))^+=G^+\} $.
Similarly,
we write $z^*(G)$ and $z^*(G^+)$ for the cut $z(G^+)$ build with
respect to $(M,<,\cdot )$.

\medskip
\medskip

Finally, let again  $T$ be a polynomially bounded o-minimal
expansion of fields with archimedean prime model, such that
$T$ has quantifier elimination and a universal system of axioms.
Let $K\subseteq
\R $ be the field of exponents of $T$ (recall that $K$ is the set
of all $\lambda \in \R $ such that the power function $x\mapsto
x^\lambda $, $x>0$ is definable in all models of $T$; so $K=\Q $
if $T$ is the theory of real closed fields). Let $p$ be a cut of a
model $M$ of $T$. Then every subset of $M^n$, $0$-definable in the
structure $(M,p^L)$ is quantifier free definable without
parameters in the language obtained from $\SL $ by introducing
names for the following subsets of $M$ and $M^2$ respectively:
\medskip
\beg {(-------------------)}
\item {}$p^L, G(p)$ and $V(p)$, all contained in $M$.
\smallskip
\item {}$\{ (a,b)\in M^2\st a+b\mal (V(p)^+)^\eta <p\}$ for $\eta \in \{ -1,1\} $.
\smallskip
\item {}$\{ (a,b)\in M^2\st a+b\mal (V(p)^+)^\eta <z(p)\}$ for $\eta \in \{ -1,1\} $.
\smallskip
\item {}$\{ (a,b)\in M^2\st a+b\mal (G(p)^+)^\lambda <p\}$ for $\lambda \in K$.
\smallskip
\item {}$\{ (a,b)\in M^2\st a+b\mal (G(p)^+)^\lambda <z(p)\}$ for $\lambda \in K$.
\smallskip
\item {}$\{ (a,b)\in M^2\st a+b\mal z^*(G(p))^\lambda <p\}$ for $\lambda \in K$.
\smallskip
\item {}$\{ (a,b)\in M^2\st a+b\mal z^*(G(p))^\lambda <z(p)\}$ for $\lambda \in K$.
\smallskip
\endbeg
\medskip
In other words, $Th(M,p^L)$ has quantifier elimination
if we add names for these sets to the language $\SL (\SD )$.\medskip

In section \ref {section:examples} we show that,\smallskip
\beg {(a)}
\item{(a)}there is a model $M$ of $T$ and a cut $p$ of $M$,
such that $Th(M,p^L)$ does not have quantifier elimination if we only
add names for each subset of $M$, definable in $(M,p^L)$
(example \ref {optimality}).\smallskip
\item{(b)}there is a model $M$ of $T$ and a cut $p$ of $M$, such that
$Th(M,p^L)$ does not have quantifier elimination in any language
containing $\SL (\SD )$ enlarged by only finitely many symbols
(example \ref {non finite}).
\endbeg
\medskip

If $T$ is not polynomially bounded and $p$ is the upper edge of a
convex subgroup of $M$ we still can determine the theory of
$(M,p^L)$ and obtain a quantifier elimination result; provided
that the invariance ring of $G$ is definable in an expansion
$(M,V)$  for some $T$-convex valuation ring $V$ of $M$. This will
be our first task in section \ref{section:groups}

\title{}{*}{Heirs of Cuts and Modelcompleteness Results.}{}
\label{section:review}{\the \paragraphcounter}

In this section we recall notions and results needed from \zitat{[Tr]}.\einzug
$T$ always denotes a complete, o-minimal expansion
of a field in the language $\SL $.  If $M\prec N$ are models of
$T$ and $A\subseteq N$, then we write $M\< A\> $
for the definable closure of $M\cup A$ in $N$ (which itself is
an elementary restriction of $N$).\einzug
If $f:M^n\lra M$
is a definable map of a model $M$ of $T$ and
$p\in S_n(M)$ is an $n$-type, then $f$ extends to a map
$S_n(M)\lra S_1(M)$, which we denote by $f$ again.
By o-minimality, the set $S_1(M)$ of 1-types of $M$
can be viewed as the disjoint union
of $M$ with the cuts of $M$.\einzug
If $p$ and $q$ are cuts of $M$,
then we write $p\sim q$ if there is a definable map
$f:M\lra M$ with $f(p)=q$. By \zitat{[Ma]}, Lemma 3.1, the relation $\sim $
is an equivalence relation between the cuts of $M$.
\medskip

For a certain class of cuts the elementary theory over a set of parameters can easily be described:\medskip

\Definition {}
A cut $p$ of a model $M$ of $T$ is called \notion {dense} if
$p$ is not definable and $M$ is dense in $M\< \alpha \> $,
for some realization $\alpha $ of $p$.
In \ref {characterization of dense types} other descriptions
of density are given.
Important for us is: $p$ is dense if and only if $p$ is not definable
and $G(p)=0$ (these cuts are also called Veronese cuts).\medskip

\fact {Theorem}{qe dense II}{ Let $A\prec M,N$ be models of $T$ and
let $p,q$ be dense cuts of $M,N$ respectively. Then $(M,p^L)\equiv
_A(N,q^L)$ if and only if $p\uhr A=q\uhr A$. Hence if $T$ has
quantifier elimination and a universal system of axioms, then
the $\SL (\SD )$-theory $T^{dense}$ which expands $T$ and says
that $\SD $ is a set $p^L$ of a dense cut $p$, has quantifier elimination. }
\pr by \ref {qe dense}. Density can be axiomatized in the language
$\SL (\SD )$, since we can say that the invariance group of the cut $p$ is trivial.
Special cases of this theorem can also be found in \zitat{[MMS]}.\ende
\medskip

\Definition {def gheir}
Let $M,N$ be models of $T$ and let $A\prec M,N$ be a common
elementary substructure.
Let $p_1,...,p_n$ be mutually distinct
cuts of $M$ and let $q_1,...,q_n$ be mutually distinct cuts of $N$.
Let $\CD _1,...,\CD _n$ be new unary
predicates. We say the tuple
$(q_1,...,q_n)$ is an \notion {heir} of $(p_1,...,p_n)$ over $A$ if the
following condition holds:
if $\phi (x_1,...,x_k)$ is an $\SL $-formula with parameters from $A$,
$\psi (x_1,...,x_k)$ is a quantifier free
$\SL (\CD _1,...,\CD _n)$-formula with parameters from $A$ and if\\
\zl{$(N,q_1^L,...,q_n^L)\models \exists \\ x\ \phi (\\ x)\land \psi
(\\ x)$}
then\\
\zl{$(M,p_1^L,...,p_n^L)\models \exists \\ x\ \phi (\\ x)\land \psi
(\\ x)$}\medskip
If $A=M$ we say that $(q_1,...,q_n)$ is an heir of $(p_1,...,p_n)$.
\medskip\medskip

\fact {Theorem}{charact of definable in Tconvex II}{
Let $M$ be a model of $T$, let $V$ be a $T$-convex valuation ring of
$M$ and let $A$ be a subset of $M$.
If $p$ is a nondefinable cut of $M$,
such that $p^L$ is A-definable in
$(M,V)$, then there is an $A$-definable map $f:M\lra M$ and elements
$a,b\in \dcl (A)$, $b\neq 0$
such that $f(V^+)=V(p)^+$ and such that
$p=a+b\mal V(p)^+$ or $p=a+b\mal \Dm (p)^+$.
Moreover $V\subseteq G(p)$ or $G(p)\subseteq \Dm (V)$,
where $\Dm (V)$ denotes the maximal ideal of $V$.
}\pr
By \ref {charact of definable in Tconvex}.\ende

\fact {Definition}{}{}
We say that $p$ has the \notion {signature alternative} if $\sign p\neq 0$
or if $p\not \sim \hat p$ and $q\not \sim \hat q$ where $q$ denotes the
unique extension of $p$ on $M\< \alpha \> $ ($\alpha \models \hat p$).
(c.f. \ref {signature alternative}).\medskip

\fact {Theorem}{signature alternatives II}{
Let $p$ be a cut of a model of $T$, such that $V(p)^+\sim V^+$ for some
$T$-convex valuation ring $V$ of $M$.
\beg {(ii)}
\ro The signature alternative holds for $p$. In particular
$\sign p=0$ implies $p\not \sim \hat p$.
\ro The signature alternative holds for $\a p\a $ with respect to $M^{>0}$.
In particular $\sign ^*p=0$ implies $p\not \sim G^*(p)^+$.
\endbeg
}\pr
By \ref {signature alternatives}.\ende

\medskip
\fact {Theorem}{multi heir}{
Let $T$ be polynomially bounded and let $p$
be a cut of a model $M$ of $T$. Let $q$ be a cut extending
$p$ of some $N\succ M$. Then\smallskip
\beg {(iii)}
\ro If $p=a\pm \hat p$ for some $a\in M$ then $q$ is an heir of $p$ if and only if
$q=a\pm \hat q$ and $\hat q$ is an heir of $\hat p$.\smallskip
\ro If $\sign p=0$ then $q$ is an heir of $p$ if and only if
$\hat q$ extends $\hat p$.\smallskip
\ro If $p=\hat p$ and $\sign ^*p=0$, then
$q$ is an heir of $p$ if and only if $q=\hat q$ and
$V(q)$ lies over $V(p)$.
\ro If $p$ is a dense cut, i.e. $\sign p=0$ and $G(p)=\{ 0\} $,
then $(V(q)^+,\hat q,q)$ is an heir of $(V(p)^+,\hat p,p)$
if and only if $G(q)=\{ 0\} $.\smallskip
\ro If $\sign p=0$, $G(p)\neq \{ 0\} $ and $\sign ^*\hat p=0$, then
$(V(q)^+,\hat q,q)$ is an heir of $(V(p)^+,\hat p,p)$ if and only
if $\hat q$ extends $\hat p$ and $V(q)$ lies over $V(p)$.\smallskip
\ro If $\sign p=0$, $G(p)\neq \{ 0\} $ and $\sign ^*\hat p\neq 0$,
then $(V(q)^+,\hat q,q)$ is an heir of $(V(p)^+,\hat p,p)$ if and only
if $\hat q$ extends $\hat p$, $V(q)$ lies over $V(p)$
and if there is some $a\in M$ such that $\hat q$ is an edge
of $a\mal V(q)^{*,>0}$.\smallskip
\endbeg
}\pr
(i) is easy and can be found in
\ref {general reduction}(i).\einzug
For the remaining statements we need that $T$ is polynomially bounded.
This means that all convex subrings of all models
of $T$ are $T$-convex. See \zitat{[Tr]}, section 4
for a summary of $T$-convex valuation rings.
With this information we can apply
\ref {signature alternatives II} for all cuts
of all models of $T$. Therefore:\medskip
(ii) holds by \ref {general reduction}(iii).
(iii) holds by \ref {general reduction}(iii)
applied to the o-minimal structure induced by $M$
on the multiplicative group of positive elements of $M$
(c.f. \zitat{[Tr]}, (5.3)).
(iv) holds by \ref {characterization of dense types}.
(v) and (vi) hold by \ref {pphatVp}.\ende

\fact {Theorem}{mc II}{
Let $T$ be model complete and
let $f(y,x_1,...,x_n)$ be a 0-definable map.\\
Let $\valO $, $\SG $, $\SZ ,\SZ ^*$ and $\SD $ be new unary
predicates and let
$c_1,...,c_n$ be new constants with respect to $\SL $. Let
$\eps ,\delta \in
\{ -1,0,1\} $ and let $\SL ^*$
be the language $\SL (\valO ,\SG, \SZ ,\SZ ^*, \SD
,c_1,...,c_n)$. Let $T^*$ be the $\SL ^*$-theory which
extends $T$ and which says the following
things about a model\\ $(M,V,G,Z,Z^*,D,d_1,...,d_n)$:
\medskip
\beg {(a)}
\lsmall $D=p^L$ for some cut $p$ of $M$, $p$ neither dense nor definable
with $\sign p=\delta $.\smallskip
\lsmall $V=V(p)$.\smallskip
\lsmall $G=G(p)$ and $\sign ^*G^+=\eps $.\smallskip
\lsmall $f(V(p)^+,d_1,...,d_n)>0$ is the upper edge of a $T$-convex
valuation ring of $M$. This is axiomatizable by the axioms
$g(f(V(p)^+,d_1,...,d_n))\leq f(V(p)^+,d_1,...,d_n)$ for each continuous,
0-definable map $g$.\smallskip
\lsmall $Z=\{ a\in M\st a+G^+=p\text { or }a-G^+=p\}$ and\\
$Z^*=\{ a\in M^*\st a\mal V=G\text { or }a\mal \Dm =G\} $.
\endbeg\smallskip
Then\smallskip
\beg {(ii)}
\ro $T^*$ is consistent if and only if there
is a model $M$ of $T$, a convex valuation ring $V$ of $M$ and
$d_1,...,d_n\in M$ such that $f(V^+,d_1,...,d_n)$ is the upper edge of a
$T$-convex valuation ring.\smallskip
\ro If $T^*$ is consistent then
$T^*$ is model complete.
\endbeg
}\pr By \ref {mc}.\ende
\medskip

\fact {Proposition}{obtain arbitrary signs II}{
Let $T$ be an o-minimal expansion of fields, let $M$ be a model of $T$
and let $G$ be a convex subgroup of $(M,+,<)$.
Let $V$ be a convex valuation  ring of $M$ with $V\subseteq V(G)$.
Let $\SC $ be a set of cuts of $M$ with $G\subseteq G(p)$
for all $p\in \SC $. Let $\eps ,\delta \in \{ -1,0,1\} $.
\smallskip
\einzug
Then there are an elementary extension $N$ of $M$ with
$\dim N/M=\aleph _0+\a \SC \a $ and a convex subgroup
$H$ of $(N,+,<)$, such that $H\cap M=G$, $\sign ^*H^+=\eps $,
$V(H)^+$ is an heir of $V^+$ and
such that for each $p\in \SC $ there is a cut $q$ of $N$
extending $p$ with $\sign q=\delta $ and $G(q)=H$.
\smallskip
\einzug
Moreover, if $\eps =0$, then we can choose $H$ so that in addition,
$V(H)$ is the convex hull of $V$ in $N$.
}
\pr By \ref {obtain both signature 0} and \ref {obtain arbitrary signs}.\ende

\title{}{*}{The elementary theory of convex subgroups}{}
\label{section:groups}{\the \paragraphcounter}

We start with a model theoretic preparation:\medskip

\fact {Lemma}{preliminary glue}{
Let $\SL \subseteq \SL ^*$ be languages and let $M^*,N^*$ be $\SL
^*$-structures. Let $A$ be a common $\SL $-substructure of
$M:=M^*\uhr \SL $ and
$N:=N^*\uhr \SL $, such that $Th (M,A)=Th (N,A)$.
Then there are $\SL ^*$-structures $\tilde M,\tilde N$ with the same
universe, such that $\tilde M\uhr \SL =\tilde N\uhr \SL $
together with elementary $\SL ^*$-embeddings $\phi :M^*\lra \tilde M$,
$\psi :N^*\lra \tilde N$, with $\phi (a)=\psi (a)$ for all $a\in A$.
}\pr
Since $Th (M,A)=Th (N,A)$, there is an $\SL $-structure
$K$ and elementary $\SL $-embeddings $\phi _0:M\lra K$,
$\psi _0:N\lra K$ such that $\phi (a)=\psi (a)$ for all $a\in A$.
We may assume that $\phi $ and $\psi $ are inclusions, hence we assume
that $M,N\prec K$. By further extending $K$ we may assume that
$K$ is $\kappa $-resplendent, where
$\kappa :=(2^{\a \SL \a +\a M\a +\a N\a })^+$
(see [Poi], 9.c.).
Since $K$ is $\kappa $-resplendent we can expand
$K$ to $\SL ^*$-structures $\tilde M$ and $\tilde N$,
such that $M^*\prec \tilde M$ and $N^*\prec \tilde N$.\ende

\fact {Theorem}{elementary equivalent for groups general}{
Let $A\prec M,N$  be models of $T$
and let $G\subseteq M,H\subseteq N$ be proper convex
subgroups (thus $G\neq 0,M$ and $H\neq 0,N$).
Suppose there is an $A$-definable map $f:A\lra A$ such that
$f(V(G)^+)$ and $f(V(H)^+)$ are the upper edges of a $T$-convex
valuation ring of $M,N$ respectively.
Then
$$(M,G)\equiv _A(N,H)$$\medskip
if and only if
the following three conditions hold:\medskip
\beg {(a)}
\lsmall $V(H)\cap A=V(G)\cap A=:V_0$, $Z^*(G)\cap A=Z^*(H)\cap A$,
$H\cap A=G\cap A=:G_0$,\smallskip
\lsmall $\sign ^*G^+=\sign ^*H^+$ and\smallskip
\lsmall If $V_0^+\sim G_0^+$, then there is some $A$-definable map
$g:A\lra A$ with
$g(V_0^+)=G_0^+$ such that\\
\zl{$g(V(G)^+)<G^+\iff g(V(H)^+)<H^+$.}
\endbeg\medskip
}\pr
Since $(M,V(G))\equiv _A(N,V(H))$ we can apply
\ref {preliminary glue} .
Therefore we may assume that $M=N$ and $V(H)=V(G)=:V$.
Let $W$ be the $T$-convex valuation ring of $M$ with upper edge
$f(V^+)$. By redefining $f$ outside
an $A$-definable interval containing $V^+$ in a suitable way,
we may assume that $f$ is an $A$-definable
homeomorphism.
Furthermore we may assume that $(M,H,G)$ is $\a A\a ^+$-saturated
and $\a A\a ^+$-strong homogeneous.\\
\case 1
$Z^*(G)\cap A=Z^*(H)\cap A\neq \emptyset $.
Let $a\in Z^*(G)\cap A$ and say $\sign ^*G^+=\sign ^*H^+=1$.
Then $G=a\mal V$ and $H=a\mal V$, thus
$G=H$.\\
\case 2
$Z^*(G)\cap A=Z^*(H)\cap A=\emptyset $.\einzug
Suppose first that $A\subseteq W$. By
\ref {charact of definable in Tconvex II} we have $W\subseteq V$.
Since $f^{-1}(W^+)=V^+$ it follows that
$f^{-1}$ is a strictly increasing $A$-definable homeomorphism, so
$f^{-1}(W^+)=W^+$. This shows that $V=W$ is $T$-convex.
From $A\subseteq V$ we get $A=V_0=V(G_0)$.
Since $G,H\not \in \{ V,\Dm (V)\} $ and $(M,V)$ is $\a A\a ^+$-saturated,
there is some $\alpha \in M$ such that
$A\< \alpha \> \not \subseteq V$ and $G\cap A\< \alpha \> =H\cap
A\< \alpha \> \in \{ 0,A\< \alpha \> \} $. In particular
$Z^*(G)\cap A\< \alpha \> =Z^*(H)\cap A\< \alpha \> =\emptyset $ and
conditions $(a),(b)$ hold for $A\< \alpha \> $ instead of $A$.
Since $G\cap A\< \alpha \>  $ is definable and
$V\cap A\< \alpha \> $ is a proper convex valuation ring of $A\< \alpha
\> $ we have $V\cap A\< \alpha \> ^+\not \sim G_0^+$
and condition (c) holds, too.
Hence it is enough to prove the theorem under the assumption of case 2
and the assumption $W\cap A\proper A$; note that
$V_0$ is also proper, since $f(G^+)=V_0^+$.\einzug
By \ref {charact of definable in Tconvex II} and our assumption in case
2 we have $G^+,H^+\not \sim _AV^+$.\\
\case {2.1}
$\sign ^*G^+=\sign ^*H^+\neq 0$.\\

\claim There are $m,m'\in M$ with $G=m\mal V$, $H=m'\mal V$
such that $t^{(M,V)}(m/A)=t^{(M,V)}(m'/A)$.\\
\case {2.1.1}
$V_0^+\not \sim G_0^+$.
By saturation there are $m,m'\in G$ with $G=m\mal V$, $H=m'\mal V$ such
that $m$ and $m'$ realize $G_0^+$. Since $V_0^+\not \sim G_0^+$, the cut
$G_0$ is not $A$-definable in $(A,V_0)$
(by \ref {explicit description}).
Hence there is a unique 1-type
$\xi $ of $(A,V_0)$ which extends the cut $G_0^+$.
Since $W\cap A\proper A$ we have $(A,V_0)\prec (M,V)$. Thus
$t^{(M,V)}(m/A)=\xi =t^{(M,V)}(m'/A)$.\\
\case {2.1.2} $V_0^+\sim G_0^+$.
By item (c) there is an $A$-definable map $g:M\lra M$ such that
$g(V_0^+)=G_0^+$ and such that $g(V^+)<G^+\iff g(V^+)<H^+$.
Now there are exactly two 1-types of $(A,V_0)$ extending the cut
$G_0^+$. One contains the statement $x<b\mal g(V_0^+)^\eps $,
the other one contains the statement $x>b\mal g(V_0^+)^\eps $.
So in order to prove the claim it is enough to find realizations $m,m'\in M$
of $G_0^+$ with $G=m\mal V,H=m'\mal V$ such that $g(V^+)<m$ iff
$g(V^+)<m$. We may assume that $G\subseteq H$.
As $G^+,H^+\not \sim _AV^+$
we get $g(V^+)<G^+$ or $H^+<g(V^+)$.
If $g(V^+)<G^+$ we take $m,m'\in M$ with $g(V^+)<m,m'$ such that
$G=mV$ and $H=m'V$.
If $H^+<g(V^+)$, then by saturation
there are $m\in
G,m'\in H$ such that $m,m'$ realize $G_0^+$ and $G=m\mal V$, $H=m'\mal
V$. In any case $m,m'$ are the elements we are
looking for and the claim is proved.\einzug
Now the theorem in the case 2.1. follows since $(M,V)$ is
$\a A\a ^+$-strong homogeneous: there is an $A$-automorphism
$\sigma :(M,V)\lra (M,V)$ with $\sigma (m')=m$.
Hence $\sigma (H)=G$ and $(M,G)\cong _A(M,H)$.\\
\case {2.2}
$\sign ^*G^+=\sign ^*H^+=0$.\\
Let $\alpha ,\ \beta $ be realizations of $G^+,\ H^+$ respectively.
Let $G_1,H_1$ be the largest convex subgroups of
$M\< \alpha \> $, $M\< \beta \> $ lying over $G$, $H$
respectively.
Since $\sign ^*G^+=0$ we know that $\sign ^*G_1^+=1$
and $Z^*(G_1^+)$ is the
set of realizations of $G^+$ in $M\< \alpha \> $. Since $\sign ^*H^+=0$
we have that $\sign ^*H_1^+=1$ and $Z^*(H_1^+)$ is the set of
realizations of $H^+$ in $M\< \beta \> $.
Moreover the requirement on the map $f$ remain valid.
Item (c) with the same map $g$ remains valid, since by
\ref {signature alternatives II}(ii) we have $V^+\not \sim G^+,H^+$.
So we can apply case 2.1 and
get $(M\< \alpha \> ,G_1)\equiv _A(M\< \beta \> ,H_1)$, hence
$(M\< \alpha \> ,V(G_1),G_1)\equiv _A(M\< \beta \> ,V(H_1),H_1)$.
Because $\sign ^*H^+=0$ we know that $(M,V,H)$ is existentially closed
in $(M\< \beta \> ,V(H_1),H_1)$ from \ref {multi heir}.
From these assertions we get
$(V,G)\equiv \> _{\exists ,A}(V,H)$. Since $T^0$ is model complete
relative $\SL $, it follows that $(M,V,G)\equiv _A(M,V,H)$, in
particular $(M,G)\equiv _A(M,H)$.\ende

\fact {Corollary}{elementary equivalent for groups Tconvex}{
Let $A\prec M,N$  be models of $T$
and let $G\subseteq M,H\subseteq N$ be proper convex
subgroups.  Suppose the following conditions hold:\medskip
\beg {(a)}
\lsmall $V(H)\cap A=V(G)\cap A=:V_0$, $Z^*(G)\cap A=Z^*(H)\cap A$,
$H\cap A=G\cap A=:G_0$,\smallskip
\lsmall $\sign ^*G^+=\sign ^*H^+$,\smallskip
\lsmall $V(G)$ and $V(H)$ are $T$-convex.\smallskip
\lsmall $V(G_0)$ is $T$-convex.\smallskip
\endbeg
Then
$$(M,G)\equiv _A(N,H).$$
}\pr
We only have to prove
condition (c) of \ref {elementary equivalent for groups general}.
We may again assume that $N=M$ and $V:=V(G)=V(H)$ by
\ref {preliminary glue}.
So assume that $V_0^+\sim G_0^+$.
Since $V_0$ and $V(G_0)$  are $T$-convex, theorem
\ref {charact of definable in Tconvex II} implies that
$V_0=V(G_0)$ and the group $G_0$ is of the form
$b\mal V_0^\eps $ for some $b\in A$ and some $\eps \in \{ \pm 1\} $.
We claim that the map $g(x):=b\mal x^\eps $ fulfills
$g(V^+)<G^+\iff g(V^+)<H^+$ as required in
\ref {elementary equivalent for groups general}(c).
If $g(V^+)=G^+$, then $Z^*(G)\cap A\neq \emptyset $ by
\ref {charact of definable in Tconvex II}, hence $H=G$ by items (a) and
(b). So we may assume that $g(V^+)\neq G^+$ and similar $g(V^+)\neq H^+$.
Suppose $G^+<g(V^+)<H^+$, hence
${1\over b}G^+<(V^+)^\eps <{1\over b}H^+$.
Since $V(G)=V(H)=V$ this is only possible if
${1\over b}G^+<1<{1\over b}H^+$ in contradiction to $b\in A$ and
$G\cap A=H\cap A$.\ende

\factname \ref {elementary equivalent for groups Tconvex}
does not give a quantifier elimination result directly.
The reason is that
\ref {elementary equivalent for groups Tconvex}(d) is not an elementary
statement in the structure $(M,G)$.
Here is an example, which shows that we can not drop assumption
(d) in \ref {elementary equivalent for groups Tconvex}.
Let $T$ be the theory of $\R _{\exp }$ and let $\omega >\R $.
Let $A:=\R \< \omega \> $ and let $V_0$ be the convex hull
of $\R $ in $A$. Let $G_0$ be the convex valuation ring
of $A$ with upper edge $\exp (\omega \mal V_0^+)$.
Let $(M,V)$ be an elementary extension of $(A,V_0)$
such that there are realizations $\alpha ,\beta \in M$ of $G_0^+$
with $\alpha <\exp (\omega \mal V^+)<\beta $.
Finally let $G:=\alpha \mal V$ and $H:=\beta \mal V$.
Then $(M,G)\not \equiv _A(M,H)$, items (a) and (b) of
\ref {elementary equivalent for groups Tconvex} holds and
$V(H)=V(G)=V$ is $T$-convex.
By \ref {obtain arbitrary signs II} it is also possible to get a counter example
from this situation in the case $\sign ^*G^+=\sign ^*H^+=0$.\medskip

If each convex valuation ring of a model of $T$ is $T$-convex, we
don't need condition (c) and (d) of
\ref {elementary equivalent for groups Tconvex}.
These theories are precisely the polynomially bounded theories with
archimedean prime model. We briefly state the definitions and known
facts which we use.\medskip

\medskip
\Definition {def poly bound}
An o-minimal expansion of a real closed field is called
\notion{polynomially bounded} if for each 0-definable map
$f$ there is some $n\in \N $ with $f(x)\leq x^n$ for sufficiently large
$x$.\medskip
\medskip
\fact {Proposition}{}{
If $T$ is an o-minimal expansion of ordered fields, then
the following are equivalent:\medskip
\beg {(ii)}
\ro $T$ is polynomially bounded and has an archimedean prime model.\smallskip
\ro Each convex valuation ring of a model of $T$ is $T$-convex.\smallskip
\endbeg
}\pr
This follows easily from \zitat {[vdD-L], (4.2)}.\ende

\medskip
{\bf For the rest of this paper we work with
a polynomially bounded $\SL $-theory $T$
which has an archimedean prime model.}\medskip

Hence from now on every convex valuation ring of a
model of $T$ is $T$-convex. This fact and
\ref {strong alternative} below gives us enough
information to determine the theory of all
$(M,p^L)$ in the language $\SL (\SD )$ over some parameter set $A$,
where $p$ runs through the cuts of $M$.\medskip

\medskip
\fact {Corollary}{elementary equivalent for groups}{
Let $A\prec M,N$  be models of $T$
and let $G\subseteq M,H\subseteq N$ be proper convex
subgroups. Then
$$(M,G)\equiv _A(N,H)$$ if and only if\medskip
\beg {(a)}
\lsmall $V(H)\cap A=V(G)\cap A=:V_0$, $Z^*(G)\cap A=Z^*(H)\cap A$,
$H\cap A=G\cap A=:G_0$ and\smallskip
\lsmall $\sign ^*G^+=\sign ^*H^+$.
\endbeg
}
\pr By \ref {elementary equivalent for groups Tconvex}.\ende
\medskip
It is not difficult to see that convex subgroups of $(M,+,<)$ are
in 1-1 correspondence with pairs $(V,\xi )$, where $V$ is a convex
valuation ring of $M$ and $\xi $ is a cut of $\Gamma _V$ with
$G(\xi )\in \{ 0,\Gamma _V\} $. If $G$ corresponds to $(V,\xi )$
and $v:M\lra \Gamma $ denotes the valuation corresponding to $V$
then \factname \ \ref {elementary equivalent for groups} says that
the theory of $(M,G)$ with parameters from $A$ is determined by
the theory of $(\Gamma ,\xi ^L)$ with parameters in $v(A^*)$. For
the theory of cuts coming from residue fields of convex valuation
rings see \ref {from res field} below. We do not make use of this
point of view.
\medskip
\medskip

\fact {Corollary}{qe groups}{ Let $\valO ,\SZ ^*$ and $\SG $ be
new unary predicates. Let\\  $\SL ^{group}:=\SL (\valO ,\SZ ^*,\SG
)$. For $\eps \in \{ -1,0,+1\} $ let $T^\eps $ be the $\SL
^{group}$-theory which extends $T$ and says the following things
about a model $(M,V,Z^*,G)$ of $T^\eps $:\medskip
\beg {(ii)} \ro $G$ is a
convex subgroup of $M$ with $\sign ^*G^+=\eps $.\smallskip
\ro $V=V(G)$ and
$Z^*=Z^*(G)$.
\endbeg\medskip
\einzug
Suppose $T$ has quantifier elimination and a universal system of axioms.
Then $T^\eps $ has quantifier elimination.
}
\pr By \ref {elementary equivalent for groups}.\ende

\fact {Corollary}{1-types of (M,G) over A}{ Let $G\subseteq M$ be
a convex subgroup and let $A\prec M$. If $p$ is a cut of $M$, such
that $p^L$ is definable in $(M,G)$ with parameters from $A$ and
$\alpha \in M$, then the 1-type $t^{(M,G)}(\alpha /A)$ is uniquely
determined by the following data:
\beg {(a)} \lsmall The cut
$t(\alpha /A)$. \lsmall $\alpha \in f(r)^L$, where $f:M\lra M$ is
$A$-definable, $r\in \{ V(G)^+,z^*(G^+),G^+\} $\\ and $f(r\uhr
A)=t(\alpha /A)$. \lsmall $\alpha \in f(r)^R$, where $f:M\lra M$
is $A$-definable, $r\in \{ V(G)^+,z^*(G^+),G^+\} $\\ and $f(r\uhr
A)=t(\alpha /A)$.
\endbeg
}\pr
By \ref {qe groups}.\ende

\title{}{*}{The elementary theory of Cuts.}{}
\label{section:cuts}{\the \paragraphcounter}

Again $T$ denotes a polynomially bounded o-minimal expansion of
fields with archimedean prime model. First we recall
the valuation property and reformulate it for our purposes.
\medskip\einzug
Let $\SP $ denote the prime model of $T$. Since $\SP $ is assumed
to be archimedean we may assume that the underlying field of $\SP
$ is an ordered subfield of $\R $. Hence the field $\R $ can be
extended to an elementary extension $\SR $ of $\SP $.\einzug A
\notion{power function} is a definable endomorphism $f$ of the
multiplicative group $\SP ^{>0}$. The extension of $f$ to $\SR $
is of the form $x^\lambda $ for some exponent  $\lambda \in \R $
and $f$ is uniquely determined by $\lambda $. The set $K$ of all
these exponents $\lambda $ is an ordered subfield of $\R $. If $V$
is a convex valuation ring of a model $M$ of $T$, then the value
group $\Gamma $ of $V$ is a $K$-vector space with multiplication
$\lambda\mal v(m):=v(m^\lambda )$ for $m\in M,\ m>0$.
All this is from
\zitat {[Mi1]}; an explanation can also be found in
\zitat {[vdD2], \S 3}.\medskip

\medskip

\fact {Theorem}{valuation property}{ (Valuation Property)\\
Let $M\prec N$ be models of $T$
with $\dim N/M=1$. Let $V\subseteq N$ be $T$-convex and let
$v:N\lra \Gamma \cup \{ +\infty \} $ be the valuation according to $V$.
Suppose $v(M)\neq v(N)$. Then\medskip
\beg {(a)}
\lsmall If $\alpha \in N\setminus M$ then there is some $a\in M$ such
that $v(\alpha -a)\not \in v(M)$.\smallskip
\lsmall If $\gamma \in v(N)\setminus v(M)$ then $v(N^*)=v(M^*)\oplus
\gamma \mal K$.
\endbeg
}\pr
(a) is \zitat{[vdD-S],9.2} and (b) is \zitat{[vdD2],5.4}.\ende

We make use of another formulation of
\ref {valuation property}:\medskip

\fact {Proposition}{move group}{
Let $G\subseteq M\models T$ be a convex
subgroup and
let $f:M\lra M$ be $M$-definable and non constant in each neighborhood
of $G^+$.
Then there is a unique exponent $\lambda $ with the following property:
there are $a,b\in M$ such that
for all $N\succ M$ and all convex subgroups $H$ of $N$ with $H\cap M=G$
we have
$$f(H^+)=a+b\mal (H^+)^\lambda .$$
We have $\lambda \neq 0$ and in the case
$\sign ^*G^+=0$, $\lambda $  is the unique exponent with the property
$f(G^+)=c+d\mal (G^+)^\lambda $ for some $c,d\in M$.
}\pr
First we prove that there are $a,b\in M$ and an exponent $\lambda \neq
0$ as claimed.
Let $\alpha $ be a realization of $G^+$ and let
$V$ be the convex hull of $\Q $ in $M\< \alpha \> $.
If $a\in M, a>0$ and $v(\alpha )=v(a)$, then there are $n,m\in \N $
such that $a<n\alpha $ and $\alpha <ma$. Since $\alpha \models G^+$
and $G$ is a convex subgroup, this is not possible.\einzug
Therefore $v(\alpha )\not \in v(M)$ and we can apply
\ref {valuation property}.
Since $f(\alpha )\not \in M$ there is some
$a\in M$ such that $v(f(\alpha )-a)\not \in \Gamma _{V\cap M}$.
Hence there are $b\in M$ and an exponent $\lambda \neq 0$, such that
$v(f(\alpha )-a)=v(b)+\lambda v(\alpha )$.
We claim that $a,b$ and $\lambda $ are the elements we are looking
for.\einzug
To see this take $N\succ M$ and a convex subgroup $H$ of $N$ with
$H\cap M=G$. Let $\beta $ be a realization of $H^+$
and let $W$ be the convex hull of $\Q $ in $N\< \beta \> $.
Furthermore let $V:=W\cap M\< \beta \> $.
By what we have shown before we have $v(f(\beta )-a)=
v(b)+\lambda v(\beta )$. Hence $w(f(\beta )-a)=
w(b)+\lambda w(\beta )$ as well.
Since $\lambda \neq 0$, we know that $f(\beta )-a$ realizes
the same cut of $N$ as $b\mal \beta ^\lambda $.
Hence $t(f(\beta )/N)=t(a+b\mal \beta ^\lambda /N)$ and $a,b$ and
$\lambda $ have the required properties.\einzug
Now we prove the second uniqueness statement.
So assume $\sign ^*G^+=0$, $c,d\in M$ and $\mu $ is an exponent such
that $f(G^+)=c+d\mal (G^+)^\mu $.
Hence $b\mal (G^+)^\lambda =c-a+d\mal (G^+)^\mu $, thus
$G^+=({d\over b})^{1\over \lambda }\mal (G^+)^{\mu \over \lambda }$.
We write $e:=({d\over b})^{1\over \lambda }$ and
$\eta :={\mu \over \lambda }$.
Let $\alpha $ be a realization of $G^+$ and let $V$ be the convex hull
of $V(G)$ in $M\< \alpha \> $ with corresponding valuation $v$.
Since $\sign ^*G^+=0$ there is some $m\in G$ with ${\alpha \over
e\mal \alpha ^\eta }<m$. Hence $v(\alpha )=v(e\alpha ^\eta )=v(e)+\eta
\mal v(\alpha )$.
Now $\eta \neq 1$ implies $v(\alpha )={1\over 1-\eta }v(a)\in
v(M)$ in contradiction to $\sign ^*G^+=0$. So $\eta =1$, i.e. $\mu
=\lambda $, which proves the second uniqueness statement.\einzug
The first uniqueness statement concerning $\lambda $ follows from
the second one, since by
\ref {obtain arbitrary signs II} there is some
$N\succ M$ and a convex subgroup $H$ of $N$, lying over $G$
with $\sign ^*H^+=0$.\ende

Proposition \ref {move group} applied to $N=M$ gives a strengthening
of the signature alternative
(\ref {signature alternatives II}(i))
for polynomially bounded theories with
archimedean prime model:\medskip

\fact {Corollary}{strong alternative}{
Let $G$ be a convex subgroup of a model $M$ of $T$.
If $p$ is a cut of $M$ with $p\sim G^+$, then $\sign p\neq 0$.
}\ende

The multiplicative signature alternative \ref {signature alternatives II}(ii)
can not be strengthened in this way. To see an example, let
$R$ be a real closed field,
let $t$ be infinitesimal, positive over $R$ and let $S$ be the real
closure of
$R(t)$.
Let $v$ be the valuation of $S$ according to the convex hull of
$R$ in $S$ and let
$G:=\{ a\in S\ \vert \ v(a)>\sqrt 2\}
$. Let $p:=G^+$ and $q:=1+p$. Since $1+G$ is a convex subgroup of
$(S^{>0},\cdot ,<)$ we have $\sign ^*q=1$.
But $\sign ^*p=0$ and $p\sim q$.\einzug
If $T$ is an expansion of $Th (\R ,\exp )$ and $R$ is a model of $T$,
then we can apply the logarithm of $R$ to the example above. This shows
that \ref {strong alternative} fails without the assumption that
$T$ is polynomially bounded. Note that $(\R ,\exp)$ is o-minimal by \zitat {[Wi]} and that an o-minimal theory $T$ of fields is an expansion
of $(\R ,\exp)$ if $T$ is not polynomially bounded (by \zitat{[Mi2]}).
\medskip

We come to our main result:\medskip

\fact {Theorem}{general main}{ Let $T$ be polynomially bounded
with archimedean prime model. Let $M,N$ be models of $T$, $A\prec
M,N$ and let $p$, $q$ be cuts of $M,N$ respectively. Then
$$(M,p^L)\equiv _A(N,q^L)$$ if and only if the following four
conditions are fulfilled:\medskip
\beg {(iii)}
\ro $p\uhr A=q\uhr
A=:p_0$, $Z(p)\cap A=Z(q)\cap A$, $G(p)\cap A=G(q)\cap A=:G_0$,
$Z^*(\hat p)\cap A=Z^*(\hat q)\cap A$  and $V(p)\cap A=V(q)\cap
A=:V_0$. \medskip
\ro $p$ is dense if and only if $q$ is dense, $p$ is
definable if and only if $q$ is definable, $\sign p=\sign q$ and
$\sign ^*\hat p=\sign ^*\hat q$.\medskip
 \ro If $p_0\sim V_0^+$ then there
are $a,b\in A$ and $\eps \in \{ \pm 1\} $ with $a+b(V_0^+)^\eps
=p_0$ such that\smallskip
\beg {(a) }
\lsmall $a+b(V(p)^+)^\eps <p\iff
a+b(V(q)^+)^\eps <q$ and\smallskip
 \lsmall $a+b(V(p)^+)^\eps <z(p)\iff
a+b(V(q)^+)^\eps <z(q)$.
\endbeg\medskip
\ro If $p_0\sim G_0^+$ then
there is some $a\in A$ with $\a p_0-a\a =\hat p_0$
such that for all $b\in A$ and each exponent $\lambda $
with $a+b(G_0^+)^\lambda =p_0$ we have\smallskip
\beg {(a) }
\lsmall $a+b\mal (\hat p)^\lambda <p\iff
a+b\mal (\hat q)^\lambda <q$,\smallskip
\lsmall $a+b\mal (\hat p)^\lambda <z(p)\iff
a+b\mal (\hat q)^\lambda <z(q)$,\smallskip
\lsmall $a+b\mal z^*(\hat p)^\lambda <p
\iff a+b\mal z^*(\hat q)^\lambda <q$ and\smallskip
\lsmall $a+b\mal z^*(\hat p)^\lambda <z(p)
\iff a+b\mal z^*(\hat q)^\lambda <z(q)$.
\endbeg
\endbeg
}\pr
Suppose first that
$(M,p^L)\equiv _A(N,q^L)$. Clearly (i) and (ii) holds.
Item (iv) holds, since $G_0^+\sim p_0$ implies that $\sign p_0\neq 0$
(by \ref {strong alternative}). Item (iii) holds since
$V_0^+\sim p_0$ implies that there are $a,b\in A$ and $\eps \in \{
+1,-1\} $ with $a+b(V_0^+)^\eps =p_0$
(by \ref {charact of definable in Tconvex II}).\einzug
Conversely suppose that (i)-(iv) holds.
If $p$ and $q$ are definable, then
$(M,p^L)\equiv _A(N,q^L)$, since $p\uhr A=q\uhr A$, $\sign p=\sign q$ and
$Z(p)\cap A=Z(q)\cap A$.
If $p$ and $q$ are dense, then $(M,p^L)\equiv _A(N,q^L)$ follows
from $p\uhr A=q\uhr A$ and \ref {qe dense}. So we may assume that
$p$ and $q$ are neither dense nor definable.
By \ref {elementary equivalent for groups} conditions (i) and (ii) imply
$(M,G(p))\equiv _A(N,G(q))$. Since
$z(q)$ and $z^*(\hat q)$ are 0-definable in $(M,q^L)$
we may assume that $M=N$, $G(p)=G(q)=:G$ and $V(p)=V(q)=:V$ by
\ref {preliminary glue}.\smallskip
\case 1
$\sign p=\sign q\neq 0$. We do the case $\sign p=\sign q=1$, the case
$\sign p=\sign q=-1$ is similar.
By (i) and (ii) we have $z(p)=z(q)$ and $z^*(\hat p)=z^*(\hat q)$.\smallskip

\case {1.1} $z(p)\uhr A=z(q)\uhr A\neq p_0$.
Hence there are $a,b\in A$ with $p=a+G^+$ and $q=b+G^+$.
Then $b-a\in G$, so $p=q$ and we are done.\smallskip
\case {1.2}$z(p)\uhr A=z(q)\uhr A=p_0$,
hence there are no $a\in A$ and no convex subgroup $H$ of $(M,+,<)$ such
that $p=a+H^+$ or $q=a+H^+$.\einzug
Let $\alpha ,\beta \in M$ such that $p=\alpha +G^+$ and $q=\beta
+G^+$.
We prove that $t^{(M,G)}(\alpha /A)=t^{(M,G)}(\beta /A)$.
For this we have to go
through conditions (a)-(c) of \ref {1-types of (M,G) over A}.
First we have $t(\alpha /A)=t(\beta /A)=p_0$ by our assumption
of case 1.2.
Let $r\in \{ V^+,z^*(G^+),G^+\} $ and let $f:M\lra M$ be $A$-definable
with $f(r\uhr A)=p_0$.
In order to check conditions (b) and (c) of
\ref {1-types of (M,G) over A} we have to show that $f(r)<\alpha $ if and
only if $f(r)<\beta $.
Since $r$ is the edge of a convex subgroup we can use
\ref {move group} to find
$a',b\in A$ and an exponent $\lambda \neq 0$
such that $f(r)=a'+b\mal r^\lambda $. Hence we may assume that
$f(x)=a'+b\mal x^\lambda $.
Suppose $\alpha <f(r)<\beta $. We may assume that $p<q$, thus $p<\beta
$ and $p<\beta -G^+=z(q)$.\smallskip
\case {1.2.1} $r=G^+$.\einzug
Take $a\in M$ such that
$\a p_0-a\a =\hat p_0$ and such that (iv)(a)-(d) holds.
Since $\alpha -a$ and $\beta -a$ realize the same edge of the
subgroup $G(p_0)$ of $A$ we have that $(f(G^+)-a)\uhr A$
is this edge. Now $(f(G^+)-a)\uhr A=a'-a+b\mal (G_0^+)^\lambda $,
so $a'-a\in b\mal G_0^\lambda $ and
$a'-a\in b\mal G^\lambda $.
It follows $a+b(G^+)^\lambda =f(r)<\beta <q$, hence
$z(p)<\alpha <a+b(G^+)^\lambda <p<z(q)$, by the choice of $a$ and (iv)(a).
This contradicts (iv)(b).\smallskip
\case {1.2.2} $r=z^*(G^+)$.\einzug
If $z^*(G^+)\uhr A\neq G_0^+$, then $z^*(G^+)\sim _AG^+$ and
we can use case 1.2.1. If $z^*(G^+)\uhr A=G_0^+$
then the same proof as in case 1.2.1. leads to a contradiction
(where we use (iv)(c) and (d) instead of (a), (b) now).\smallskip
\case {1.2.3} $r=V^+$.\einzug
Take $a_0,b_0\in A$ and $\eps \in \{ \pm 1\} $ as in (iii).
Since $r=V^+$ and $p_0=a+b\mal (V^+)^\lambda $
we may assume that $a_0=a'$ as in the proof of
case 1.2.1.
So $f(r)=a_0+b\mal V^\lambda $ and $b\mal (V_0^+)^\lambda =G(p_0)$.
Since $a_0+b_0(V_0^+)^\eps =p_0$ it follows that
$b\mal (V_0^+)^\lambda =b_0\mal (V_0^+)^\eps $, hence
$b\mal (V^+)^\lambda =b_0\mal (V^+)^\eps $, too.
This shows that we can replace
$b$ by $b_0$ and $\lambda $ by $\eps $,
hence $f(r)=a_0+b_0\mal (V^+)^\eps $.
It follows $a_0+b_0\mal (V^+)^\eps <\beta <q$, hence
$a_0+b_0\mal (V^+)^\eps <p$ by (iii)(a). So
$z(p)<\alpha <a_0+b_0\mal (V^+)^\eps <p<z(q)$ which contradicts
(iii)(b).\smallskip\einzug
This proves that
$t^{(M,G)}(\alpha /A)=t^{(M,G)}(\beta /A)$. Since $(M,G)$ is strong
$\a A\a ^+$-homogeneous there is an $A$-automorphism
$\sigma $ of $(M,G)$ which maps $\beta $ to $\alpha $.
Hence the cut $q$ is mapped to the cut $p$ under this automorphism,
which shows $(M,p^L)\equiv _A(N,q^L)$.\medskip

\case 2 $\sign p=\sign q=0$.\einzug
Then $z(p)=p$ and $z(q)=q$.
Let $\alpha ,\beta $ be realizations of $p,q$ respectively
and let $p_1,q_1$ be the largest extensions of $p$,$q$ on
$M\< \alpha \> $, $M\< \beta \> $ respectively.\smallskip

\Claim A $(M\< \alpha \> ,p_1^L)\equiv _A(M\< \beta \> ,q_1^L)$.
To see claim A we may use case 1, since $\sign p_1=\sign q_1=1$.
Hence it is enough to check conditions (i)-(iv) for $p_1$ and
$q_1$. Condition (ii) clearly holds for $p_1$ and $q_1$, since
$\hat p_1$, $\hat q_1$ are the unique extensions of $\hat p$,
$\hat q$ respectively, so the signature does not change.\einzug
Since $\sign p=0$ the cut $z^*(\hat p_1)$ is the unique extension
of $z^*(\hat p)$ on $M\< \alpha \> $. Since $p_1=\alpha +\hat p_1$
we have $z(p_1)=\alpha -\hat p_1$ is the least extension of $p$ on
$M\< \alpha \> $. Moreover $V(p_1)$ is the convex hull of $V$ on
$M\< \alpha \> $. From these data it follows that conditions
(i),(iii) and (iv) are fulfilled for $p$ and $p_1$. Similarly they
are fulfilled for $q$ and $q_1$. Hence (i)-(iv) holds for $p_1$
and $q_1$.\medskip

\Claim B $(V^+,G^+,z^*(G^+),p)\equiv \> _{\exists ,A}(V^+,G^+,z^*(G^+),q)$.
\einzug
This follows from claim A if we know that
$(V_{q_1}^+,\hat  q_1,z^*(\hat q_1),q_1)\equiv \> _{\exists ,A}
\ (V^+,G^+,z^*(G^+),q)$. As $\sign q=0$, the cuts $G^+$ and
$V^+$ are omitted
in $M\< \beta \> $, in particular
$\sign ^*\hat q_1=\sign ^*\hat q$ and $Z^*(\hat q_1)\cap M=Z^*(\hat q)$.
So by \ref {multi heir} we even know
$$(V_{q_1}^+,\hat  q_1,z^*(\hat q_1),q_1)\equiv \> _{\exists ,M}
\ (V^+,G^+,z^*(G^+),q).$$
In order to prove the theorem we use $\ref {mc II}$.
Let
$Z_M^*:=\{ a\in M^*\st a\mal V(p)=G(p)\text { or }a\mal \Dm (p)=G(p)\} $
and
$Z_N^*:=\{ a\in N^*\st a\mal V(r)=G(r)\text { or }a\mal \Dm (r)=G(r)\}
$. Let $\DM $ be the $\SL ^*$-structure $(M,V,G,Z_M^*,p^L)$
and let $\DN $ be the $\SL ^*$-structure
$(M,V,G,Z_N^*,q^L)$. Let $\eps :=\sign ^*\hat p$.
We know that $\DN ,\DM \models T^\eps _0$.
Claim B says that for each existential
$\SL ^*(A)$-sentence $\phi $ we have $\DM \models \phi \Rightarrow
\DN \models \phi $. Since $T^\eps _0$ is model complete we get
$\DM \equiv _A\DN $, hence
$(M,p^L)\equiv _A(M,q^L)$.\ende

\medskip

In example \ref {example3} below,
we construct a situation, which shows
that (i),(ii) and (iv) of \ref {general main} do not
imply item (iii).\medskip
\medskip

\fact {Corollary}{abstract content of general main}{
There is some $a\in A$ such that
$$(M,p^L)\equiv _A(N,q^L)\iff (M,p^L)\equiv _{\{ a,b\} }(N,q^L)
\text { for each }b\in A.$$
}\pr
With the notation as in \ref {general main} we can take
$a=0$ if $\sign p_0=0$ and $a\in A$ with $\a p_0-a\a =\hat p_0$
otherwise.\ende

\medskip
By \ref {optimality} below it is possible to have
$(M,p^L)\not \equiv _A(N,q^L)$ and $(M,p^L)\equiv _{\{ b\} }(N,q^L)$
for each $b\in A$.\einzug
Special cases of \ref {general main} are formulated in the next
corollaries:\medskip
\medskip

\fact {Corollary}{sign* neq 0 with witness}{
Let $M,N$ be models of $T$, $A\prec M,N$ and let
$p$, $q$ be cuts of $M,N$ respectively.
Suppose\smallskip
\beg {(iii)}
\ro $p\uhr A=q\uhr A=:p_0$, $Z(p)\cap A=Z(q)\cap A$,
$G(p)\cap A=G(q)\cap A=:G_0$,
$Z^*(\hat p)\cap A=Z^*(\hat q)\cap A\neq \emptyset $  and
$V(p)\cap A=V(q)\cap A=:V_0$.\smallskip
\ro $p$ and $q$ are neither dense nor definable with $\sign p=\sign q$.\smallskip
\ro If $p_0\sim V_0^+$, then there are $a,b\in A$ and $\eps
\in \{ \pm 1\} $ with $a+b(V_0^+)^\eps =p_0$ such that\smallskip
\beg {(a)}
\lsmall $p<a+b(V(p)^+)^\eps \iff q<a+b(V(q)^+)^\eps $ and
\lsmall $z(p)<a+b(V(p)^+)^\eps \iff z(q)<a+b(V(q)^+)^\eps$.
\endbeg
\endbeg\medskip
Then $$(M,p^L)\equiv _A(N,q^L).$$
\medskip
}\pr
Since $\sign ^*\hat p=\sign ^*\hat q$ is implied by (i) and (ii),
we only have to check condition (iv) of \ref {general main}.
By \ref {elementary equivalent for groups} we know
$(M,G(p))\equiv _A(N,G(q))$ and
we may again assume that $M=N,\ V:=V(p)=V(q)$ and $G:=G(p)=G(q)$.
Moreover we assume that $\sign ^*G^+=1$, the case $\sign ^*G^+=-1$ is
similar.
Let $c\in Z^*(\hat p)\cap A=Z^*(\hat  q)\cap A$, so
$G=c\mal V$ and $z^*(G^+)=c\mal \Dm ^+$.\einzug
Now let $p_0\sim G_0^+$ as in \ref {general main}(iv) assumed.
Since $G_0=c\mal V_0$ it follows $p_0\sim V_0^+$, hence by (iii) there
are $a,b\in A$ and some $\eps \in \{ \pm 1\} $
with $a+b(V_0^+)^\eps =p_0$ and $p<a+b(V^+)^\eps \iff q<a+b(V^+)^\eps$.
We prove \ref {general main}(iv)(a), the other cases
are similar.
We have to show that for all $b_1\in A$ and each exponent $\lambda $
the cut
$a+b_1(G^+)^\lambda $ does not lie between $p$ and $q$.
Since $p$ and $q$ extends $p_0$,
we may assume that $a+b_1(G^+)^\lambda $ extends $p_0=a+b\mal (V^+)^\eps $,
hence $b_1(cV_0^+)^\lambda =b\mal (V_0^+)^\eps $, in other words
the $A$-definable map
$f(x):=({b_1\mal (cx)^\lambda \over b})^{1\over \eps} $ fixes $V_0^+$.
Hence $f$ fixes $V^+$, too.
Thus $b_1(cV^+)^\lambda =b\mal (V^+)^\eps $, which gives the claim by the
choice of $a$ and $b$.\ende

\smallskip
\einzug
A cut $p$ of $M$ is of the form $P^{-1}(\xi ^L)^+$
where $P:M\lra M_0\cup \{ +\infty \} $ is a real place
and $\xi $ is a cut of $M_0$ with definable invariance group if and only
if the valuation ring of $P$ is $V(p)$, $G(p)=\Dm (p)$ and
$\a p\a \leq V(p)^+$. For these cuts the theory of $(M,p^L)$ with
parameters in $A\prec M$ is determined by the theory
of $(M_0,\xi ^L)$ with parameters in $P(A)$:
\medskip

\medskip
\fact {Corollary}{from res field}{
Suppose $p$,$q$ are neither dense nor definable such that
$\sign p=\sign q$, $Z(p)\cap A=Z(q)\cap A$,
$p\uhr A=q\uhr A$, $V(p)\cap A=V(q)\cap A$,
$G(p)=\Dm (p)$, $G(q)=\Dm (q)$, $\a p\a \leq V(p)^+$ and
$\a q\a \leq V(q)^+$.
Then $(M,p^L)\equiv _A(N,q^L)$.
}\pr
Say $p,q>0$. We use \ref {sign* neq 0 with witness} and we only have
to check condition (iii) there. We do statement (a), statement (b) is
similar.
Since $0<p\leq V^+$ and $G(p)=\Dm $ we have
$\Dm _0^+\leq p_0\leq V_0^+$, so if $p_0=V_0^+$ we can
take $a=0$, $b=1$ and $\eps =1$. So we assume that $\Dm _0\leq p_0<V_0^+$.
We assume again that $M=N$ and $V=V(p)=V(q)$.
Suppose $p_0=a+b(V_0^+)^\eps $ and $p\leq a+b(V^+)^\eps \leq q$.
If $p=a+b(V^+)^\eps $ or $q=a+b(V^+)^\eps $, then $p=q$
since $Z(p)\cap A=Z(q)\cap A$ and we are done.
So we assume that $p<a+b(V^+)^\eps <q$.
Since $\Dm _0^+\leq b(V_0^+)^\eps <V_0^+$ we must have
$\Dm _0=b\mal (V_0^+)^\eps $, hence $\eps =-1$ and $b\in V_0^*$.
Since $V_0^*\subseteq V^*$ we may assume that $\a b\a =1$,
Hence $p<a\pm \Dm ^+<q$, which contradicts $G(p)=G(q)=\Dm $
and $p\uhr A=q\uhr A$.\ende
\medskip

\beginlong {
\fact {Corollary}{all signs 0}{
Let $T$ be polynomially bounded with archimedean prime model.
Let $M,N$ be models of $T$, $A\prec M,N$ and let
$p$, $q$ be cuts of $M,N$ respectively. Suppose
\beg {(iii)}
\ro $V(p)\cap A=V(q)\cap A$, $G(p)\cap A=G(q)\cap A=:G_0$
and $p\uhr A=q\uhr A=:p_0$.
\ro $\sign p=\sign q=0$ and $\sign ^*\hat p=\sign ^*\hat q=0$.
\ro If $p_0\sim V_0^+$, then there are $a,b\in A$ and some $\eps
\in \{ \pm 1\} $
with $a+b(V_0^+)^\eps =p_0$ and
$p<a+bV(p)^\eps \iff q<a+bV(q)^+$.
\ro
If $p_0 \sim G_0^+$ then
there is some $a\in A$ with $\a p_0-a\a =\hat p_0$
such that for all $b\in A$ and each exponent $\lambda $
we have $a+b\mal (\hat p)^\lambda <p\iff
a+b\mal (\hat q)^\lambda <q$.
\endbeg
Then $$(M,p^L)\equiv _A(N,q^L)$$
}\pr
By (ii) we have $z(p)=p$, $z^*(\hat p)=\hat p$,
$z(q)=q$, $z^*(\hat q)=\hat q$ and neither $p$ nor $q$ is dense.
Hence (i)-(iv) of \ref {general main} are fulfilled.\ende

\fact {Corollary}{}{
Suppose $A\subseteq V(p),V(q)$, $1+p=p,1+q=q$,
and $\sign p=\sign q=
\sign ^*\hat p=\sign ^*\hat q=0$.
Then $(M,p^L)\equiv _A(N,q^L)$ if and only
if for each exponent $\lambda >1$ we have
$\hat p^\lambda <p\iff \hat q^\lambda <q$.
Moreover, for each cut $\xi $ of exponents $>1$ , there is a model
$M\succ A$ and a cut $p$ of $M$ with
$1+p=p$, $A\subseteq V(p)$ and $\sign p=\sign ^*\hat p=0$,
such that $\lambda <\xi \iff
\hat p^\lambda <p$ ($\lambda >1$)
}\pr
Let $M_1\succ A$ and let $G_1$ be a convex subgroup of $M_1$ with
$\sign ^*G_1^+=0$ and $A\subseteq V(G_1)\subseteq G_1$.
First suppose $\xi ^L,\xi ^R$ are nonempty.
Let $\alpha \in M_2\succ M_1$ with $(G_1^+)^\lambda \leq t(\alpha
/M_1)\leq
(G_1^+)^\mu $ for all $1<\lambda <\xi <\mu $ and let
$G_2$ be the convex hull of $G_1$ in $M_2$.
Clearly $(G_2^+)^\lambda <\alpha +G_2^+<(G_2^+)^\mu $.
Finally let $M\succ M_2$ such that there is an extension
$p$ of $\alpha +G_2^+$ on $M$ with $G(p)\cap M_2=G_2$,
$\sign ^*\hat p=0$ and such that $V(p)$ lies over $V(G_2)$.\einzug
Now suppose $\lambda <\xi $ for all exponents.
Then we take $\alpha \in M_2\succ M_1$ such that $G^\lambda <\alpha $
for all $\lambda $ and proceed as above.

\underconstruction {noch fertig machen}
}\endlong

\fact {Corollary}{sign p0 = 0}{
Let $T$ be polynomially bounded with archimedean prime model.
Let $M,N$ be models of $T$, $A\prec M,N$ and let
$p$, $q$ be cuts of $M,N$ respectively. Suppose\medskip
\beg {(iii)}
\ro $p\uhr A=q\uhr A$ has signature 0,
$G(p)\cap A=G(q)\cap A$,
$Z^*(\hat p)\cap A=Z^*(\hat q)\cap A$  and
$V(p)\cap A=V(q)\cap A=:V_0$.
\smallskip
\ro $p$ is dense if and only if $q$ is dense, $p$ is definable if and
only if $q$ is definable,
$\sign p=\sign q$ and $\sign ^*\hat p=\sign ^*\hat q$.
\begincomment {
\ro $\sign p=\sign q$ and $p\uhr A=q\uhr A$ has signature 0.
\ro $(M,G(p))\equiv _A(N,G(q))$.
}\endcomment 
\endbeg\medskip
Then $$(M,p^L)\equiv _A(N,q^L)$$
}
\pr Since $p_0:=p\uhr A$ has signature 0 we have
$Z(p)\cap A=\emptyset =Z(q)\cap A$. Moreover $p_0\not \sim V_0^+$
and $p_0\not \sim G_0^+$, where $V_0=V(p)\cap A$ and $G_0=G(p)\cap A$.
Hence (i)-(iv) of \ref {general main} are fulfilled and the \factname
follows.\ende

\smallskip
Theorem \ref {general main} implies a quantifier elimination
result for the theories $T^\eps _\delta $ in an extended
language:\medskip

\fact {Corollary}{qe language}{
Suppose $T$ has quantifier elimination and a universal system of axioms.
Let $\eps ,\delta \in \{ -1,0,+1\} $.
We use the notation of \ref {mc II}; since $T$ is polynomially bounded
with archimedean prime model
the language $\SL ^*$ is $\SL (\valO ,\SG ,\SZ ,\SZ ^*,\SD )$
(no constants are needed)
and does not depend on a given 0-definable map as in \ref {mc II}.
Let $\SL ^{cut}$ be the language $\SL ^*$ together with
binary predicates $R_\eta ^1,R_\eta ^2,S_\lambda ^1,...,S_\lambda
^4$ for each exponent $\lambda \neq 0$ and $\eta \in \{ \pm 1\} $.
Let $\overline {T^\eps _\delta }$ be the $\SL ^{cut}$-theory, which extends
$T^\eps _\delta $ and which says in addition the following things about
a model $(M,p^L,G,V,...) $: for all $a,b\in M$ we have\medskip
\beg {(---------------------------------)}
\item{}$R_\eta ^1(a,b)\leftrightarrow a+b\mal (V^+)^\eta <p$,\smallskip
\item{}$R_\eta ^2(a,b)\leftrightarrow a+b\mal (V^+)^\eta <z(p)$,\smallskip
\item{}$S_\lambda ^1(a,b)\leftrightarrow a+b\mal (G ^+)^\lambda <p$,\smallskip
\item{}$S_\lambda ^2(a,b)\leftrightarrow a+b\mal (G ^+)^\lambda <z(p)$,\smallskip
$S_\lambda ^3(a,b)\leftrightarrow a+b\mal z^*(G ^+)^\lambda <p$ and
\smallskip
\item{}$S_\lambda ^1(a,b)\leftrightarrow
a+b\mal z^*(G ^+)^\lambda <z(p)$.
\endbeg\medskip

Then $\overline {T^\eps _\delta }$ has quantifier elimination.
}
\pr
By \ref {general main}, two models of
$\overline {T^\eps _\delta }$ which induce the same
$\overline \SL $-structure on a common substructure $\SA $
are elementary equivalent over $\SA $ (since
$T$ has a universal system of axiom, $\SA $ is an expansion
of a common elementary restriction $A$ of the underlying $T$-models).
This is a reformulation of quantifier elimination.\ende
\medskip
Finally, \ref {obtain arbitrary signs II}, \ref {general main} and the
corollaries above allow explicit
descriptions of the various theories $Th(M,p^L)$ where
$M$ runs through the models of $T$ and $p$ runs through the cuts
of $M$. We state here only one case:
\medskip\medskip

\fact {Corollary}{}{
Let $\xi $ be a non definable cut of the prime model of $T$.
Let $\tilde T$ be the $\SL (\SD )$-theory which extends $T$
and says the following things about a model $(M,D)$:\medskip
\beg {(3) }
\ar $D$ is the set $p^L$ of a nondefinable cut $p$ of $M$,
\ar  $p$ extends $\xi $,
\ar $p$ is not dense and $\sign p=\sign ^*\hat p=0$.
\endbeg\medskip
Then $\tilde T$ is complete.
}\pr
$\tilde T$ is consistent by \ref {obtain arbitrary signs II}.
Let $\SP $ denote the prime model of $T$ again.
If $(M,p^L)$ is a model of $\tilde T$, then
$Z(p)=Z^*(\hat p)=\emptyset $, $G(p)\cap \SP =\{ 0\}$,
$V(p)\cap \SP =\SP $ and $p\uhr \SP =\xi \not \sim O^+,\SP ^+$.\einzug
So if $(N,q^L)$ is another model of $\tilde T$, then all
conditions (i)-(iv) of \ref {general main} are fulfilled for
$A=\SP $. Hence $(M,p^L)\equiv (N,q^L)$ which shows that
$\tilde T$ is complete.\ende

\title{}{*}{Counter Examples}{}\label{section:examples}
{\the \paragraphcounter}

We give here some examples which show in which way one can not
improve \ref {general main}. $T$ is again polynomially bounded
with archimedean prime model.
\medskip

\fact {Example}{optimality}{} There are models $A\prec M,N$ of $T$
and cuts $p,q$ of $M,N$ respectively with $(M,p^L)\equiv
_b(N,q^L)$ for all $b\in A$ such that $(M,p^L)\not \equiv
_A(N,q^L)$. In other words we can have $(M,p^L)\not \equiv
_A(N,q^L)$, although for every $\SL (\SD )$-formula  $\phi (x)$ in
one variable, the sets $X:=\{ m\in M\st (M,p^L)\models \phi (m)\}
$ and $Y:=\{ n\in N\st (N,q^L)\models \phi (n)\} $ fulfill $X\cap
A=Y\cap A$. \pr We choose $A$ such that there are $a,b\in A$,
$a,b>0$ with $v_0(a)<n\mal v_0(b)<0$ for all $n\in \N $, where
$v_0$ is a convex valuation on $A$ with valuation ring $V_0$. Let
$p_0:=a+bV_0^+$ and let $M_1\succ A$ such that there are a convex
valuation ring $V_1$ lying over $V_0$ and realizations $\alpha
,\beta \in M$ of $V_0^+$ with $\alpha \in V_1<\beta $. Let
$p_1:=a+b\alpha +V^+$ and let $q_1:=a+b\beta +V^+$. Certainly $p$
and $q$ extend $p_0$. By \ref {obtain arbitrary signs II} there is some
$M\succ M_1$ and cuts $p,q$ of $M$ extending $p_1,q_1$
respectively, such that $G(p)=G(q)=$ the convex hull $V$ of $V_1$
in $M$ and such that $\sign p=\sign q=0$. So $Z(p)=Z(q)=\emptyset
$, $G(p)=G(q)=V$ and $Z^*(\hat p)=Z^*(\hat q)=(V^*)^{>0}$. In
particular conditions (i) and (ii) of \ref {sign* neq 0 with
witness} are fulfilled for any $A_0\prec A$. We have $(M,p^L)\not
\equiv _A(M,q^L)$, since $p=a+b\alpha +V^+<a+bV^+<a+b\beta
+V^+=q$.\einzug But for each $a_0\in A$ we have $(M,p^L)\equiv
_{a_0}(M,q^L)$. To see this we have to show conditions (i)-(iii)
of \ref {sign* neq 0 with witness} for $A_0:=\dcl (a_0)$. We have
already seen that conditions (i) and (ii) hold. We prove  (iii) of
\ref {sign* neq 0 with witness} for $A_0$. Suppose there are
$c,d\in A_0$ and $\eps \in \{ \pm 1\} $ such that $p\leq
c+d(V^+)^\eps \leq q$. Since $p$ and $q$ extend $p_0$ we have
$c+d(V_0^+)^\eps =a+bV_0^+$. Hence $\eps =1$ and $c+dV_0=a+bV_0$.
This means $v_0(a-c)\geq v_0(d)=v(b)<0$. Because the definable
closure of the empty set is archimedean by assumption we know that
$d$ is not 0-definable in $M$. As $\dim A_0=1$ we have $c\in \dcl
(d)$. Since $T$ is polynomially bounded there is some $n\in \N $
with $n\mal v_0(d)\leq v_0(c)$. Hence $v_0(a)<n\mal v_0(b)=n\mal
v_0(d)<v_0(c)$ and $v_0(a)=v_0(a-c)\geq v_0(d)=v_0(b)$, a
contradiction.\ende

\begincomment {
\fact {Proposition}{optimality}{
Let $T$ be polynomially bounded with archimedean prime model.
\beg {(9)}
\ar There are models $A\prec M,N$ of $T$ and cuts
$p,q$ of $M,N$ respectively with\\
$(M,p^L)\equiv _b(N,q^L)$ for all
$b\in A$ such that
$(M,p^L)\not \equiv _A(N,q^L)$.
In other words we can have
$(M,p^L)\not \equiv _A(N,q^L)$, although for
every $\SL (\SD )$-formula  $\phi (x)$ in one variable, the sets
$X:=\{ m\in M\st (M,p^L)\models \phi (m)\} $ and
$Y:=\{ n\in N\st (N,q^L)\models \phi (n)\} $ fulfill
$X\cap A=Y\cap A$.
\ar If $\SL _1\subseteq \overline \SL $ is a language extending $\SL ^*$,
which does not contain one of the $R_\eta ^k$ or $S_\lambda ^k$, then
there is a model $M$ of $T$ and a cut $p$ of $M$ such that
the $\SL _1$-theory
$Th(\DM )$ does not have quantifier elimination , where
$\DM $ denotes the $\SL (\SD )$-structure $(M,p^L)$ expanded
as in \ref {qe language}.
\ar In general, we do not have quantifier elimination of $Th(\DM )$
if $\DM $ is an expansion by definitions of $(M,p^L)$ with
finitely many new symbols.
\endbeg
}\pr
}\endcomment 

\beginlong {
A lemma in advance.\medskip

\fact {Lemma}{what means a+bV=c+dV}{
Let $R$ be a real closed field and let $V$ be a convex
valuation ring of $R$. Let $a,bc,d\in R$.
Then
$a+bV=c+dV$ if and only if $v(a-c)\geq v(b)=v(d)$.
}\pr
This certainly holds if $b\mal d=0$ or if $V=M$.
So we may assume that $bd\neq 0$ and $V$ is proper.\einzug
$\Rightarrow $.
$a+bV=c+dV$ implies $bV=dV$ - hence $v(b)=v(d)$ - and
$a-c\in bV$, i.e. $\a a-c\a <b\mal u$ with $u\in V^*$.
This implies that $v(a-c)\geq v(b)$.\einzug
$\Leftarrow $.
As $v(b)=v(d)$ we have $bV=dV$.
Since $v(a-c)\geq v(b)$, there is some $u\in V$ with
$\a a-c\a \leq b\mal u$, hence $a-c+bV=bV=dV$.\ende

\fact {Lemma}{force dim 2}{ Let $M$ be polynomially bounded with
archimedean prime model and let $V$ be a proper convex valuation
ring of $M$. Let $a,b\in V$ with $v(a)<n\mal v(b)<0$ for all $n\in
\N $. If $c,d\in M$ with $c+dV=a+bV$, then $\dim \{ c,d\} =2$.
}\pr Since $c+dV=a+bV$ we have $v(a-c)\geq v(d)=v(b)<0$. Because
the definable closure of the empty set is archimedean by
assumption we know that $d$ is not 0-definable in $M$. Suppose
$\dim \{c,d\} =1$, thus $c\in \dcl (d)$. Since $T$ is polynomially
bounded there is some $n\in \N $ with $n\mal v(d)\leq v(c)$. Hence
$v(a)<n\mal v(b)=n\mal v(d)<v(c)$ and $v(a)=v(a-c)\geq v(d)=v(b)$,
a contradiction.\ende

Now we can formulate the example:\einzug
}\endlong

\begincomment {
\fact {Example}{ex: need dim 2}{} Here is an example which shows
that we can not improve \ref {abstract content of general main}.
That is, the assumption  $(M,p^L)\equiv _b(N,q^L)$ for all $b\in
A$ does not imply $(M,p^L)\equiv _A(N,q^L)$ in general. Let $A$ be
polynomially bounded with archimedean prime model and let $a,b\in
A$, $a,b>0$ with $v_0(a)<n\mal v_0(b)<0$ for all $n\in \N $, where
$v_0$ is a convex valuation on $A$ with valuation ring $V_0$. Let
$p_0:=a+bV_0^+$ and let $M_1\succ A$ such that there is a convex
valuation ring $V_1$ lying over $V_0$ and realizations $\alpha
,\beta \in M$ of $V_0^+$ with $\alpha \in V_1<\beta $. Let
$p_1:=a+b\alpha +V^+$ and let $q_1:=a+b\beta +V^+$. Certainly $p$
and $q$ extend $p_0$. By \ref {extend to sign 0} there is some
$M\succ M_1$ and cuts $p,q$ of $M$ extending $p_1,q_1$
respectively, such that $G(p)=G(q)=$ the convex hull $V$ of $V_1$
in $M$ and such that $\sign p=\sign q=0$. So $Z(p)=Z(q)=\emptyset
$, $G(p)=G(q)=V$ and $Z^*(\hat p)=Z^*(\hat q)=(V^*)^{>0}$. In
particular conditions (i) and (ii) of \ref {sign* neq 0 with
witness} are fulfilled for any $A_0\prec A$. We have $(M,p^L)\not
\equiv _A(M,q^L)$, since $p=a+b\alpha +V^+<a+bV^+<a+b\beta
V^+=q$.\einzug But for each $a_0\in A$ we have $(M,p^L)\equiv
_{a_0}(M,q^L)$. To see this we have to show conditions (i)-(iii)
of \ref {sign* neq 0 with witness} for $A_0:=\dcl (a_0)$. We have
already seen that conditions (i) and (ii) hold. We prove  (iii) of
\ref {sign* neq 0 with witness} for $A_0$. Suppose there are
$c,d\in A_0$ and $\eps \in \{ \pm 1\} $ such that $p\leq
c+d(V^+)^\eps \leq q$. Since $p$ and $q$ extend $p_0$ we have
$c+d(V_0^+)^\eps =a+bV_0^+$. Hence $\eps =1$ and $c+dV_0=a+bV_0$.
This means $v_0(a-c)\geq v_0(d)=v(b)<0$. Because the definable
closure of the empty set is archimedean by assumption we know that
$d$ is not 0-definable in $M$. As $\dim A_0=1$ we have $c\in \dcl
(d)$. Since $T$ is polynomially bounded there is some $n\in \N $
with $n\mal v_0(d)\leq v_0(c)$. Hence $v_0(a)<n\mal v_0(b)=n\mal
v_0(d)<v_0(c)$ and $v_0(a)=v_0(a-c)\geq v_0(d)=v_0(b)$, a
contradiction.\ende

Example \ref {ex: need dim 2} shows that we can have
$(M,p^L)\not \equiv _A(N,q^L)$, although for
every $\SL (\SD )$-formula  $\phi (x)$ in one variable, the sets
$X:=\{ m\in M\st (M,p^L)\models \phi (m)\} $ and
$Y:=\{ n\in N\st (N,q^L)\models \phi (n)\} $ fulfill
$X\cap A=Y\cap A$.\medskip
}\endcomment 

\fact {Example}{non finite}{}
The following example shows that we need each exponent
in \ref {general main}. In particular
we do not have quantifier elimination of $Th(\DM )$
if $\DM $ is an expansion by definitions of $(M,p^L)$ with
finitely many new symbols.
\pr
Let $M_1\succ M_0$ such that there is a
convex subgroup $G_1$ of $M_1$ with $\sign ^*G_1^+=0$ and $M_0\subseteq
V(G_1)\subseteq G_1$. Let $\lambda \neq 0$ and let
$M_2\succ M_1$ such that
there are $b\in M_2$ with $b>M_1$ and a convex subgroup
$G_2$ of $M_2$ such that $\sign ^*G_2^+=0$, $G_2\cap M_1=G_1$ and
$V(G_2)\cap M_1=V(G_1)$. In particular
$G_2^+<b\mal (G_2^+)^\lambda $.
By compactness there is some $M_3\succ
M_2$ together with a convex subgroup $G_3$, $\sign ^*G_3^+=0$,
$G_3\cap M_2=G_2$, $V(G_3)\cap M_2=V(G_2)$ together with realizations
$\alpha ,\beta \in M_3$ of $b\mal (G_2^+)^\lambda $
such that $\alpha <b\mal (G_3^+)^\lambda <\beta $.
As $G_3^+<b\mal (G_3^+)^\lambda $ we have
$\alpha +G_3^+<b\mal (G_3^+)^\lambda <\beta +G_3^+$
and $\alpha +G_3^+,\beta +G_3^+$ are extensions of $b\mal
(G_2^+)^\lambda $.
Finally let $M\succ M_3$ such that there are extensions
$p$, $q$ of $\alpha +G_3^+,\beta +G_3^+$ on $M$
respectively, with $\sign p=\sign q=0$, $G(p)=G(q)$ lies over $G_3$,
$V(p)=V(q)$ lies over $V(G_3)$ and $\sign ^*G(p)^+=0$. Then we see:
\beg {(a)}
\lsmall Conditions (i),(ii) and (iii) of \ref {general main} hold
for $A=M_2$.
\lsmall Condition (iv) of \ref {general main} holds for $A=M_2$ and each
exponent $\mu \neq \lambda $.
\lsmall $(M,p^L)\not \equiv _{M_2}(M,q^L)$
since condition (iv) of \ref {general main} does not hold for $A=M_2$
and $\lambda $.
\lsmall $(M,p^L)\equiv _{M_0}(M,q^L)$ by (a) and (b) and
since condition (iv) of \ref {general main} holds for $A=M_0$
and $\lambda $.
\endbeg
\ende

\fact {Example}{example3}{}
The following example shows that conditions (i),(ii) and (iv) of
\ref {general main} can hold such that $(p,V(p)^+)$ is not an heir
of $(q,V(q)^+)$ over $A$ and vice versa.
Let $A\models T$ and let $V_0$ be a proper convex valuation ring of $A$.
Let $\alpha \models V_0^+$, let $W_1$ be the convex hull
of $V_0$ in $A\< \alpha \> $ and let $V_1$ be the largest convex
valuation ring
of $A\< \alpha \> $ lying over $V_0$.
Let $G_1$ be a convex subgroup of $(A\< \alpha \> ,+,<)$ such that
$V_1\subseteq V(G_1)$ and such that $G_1\subseteq W_1$ (for example
we could take $G_1=$the maximal ideal of $V_1$; if necessary we also can
choose $G_1$ such that $\sign ^*G_1^+=0$).
Furthermore let $p_1$ be a cut of $A\< \alpha \> $ such that
$W_1<p_1<V_1$ and $G(p_1)=G_1$ (as $G_1\subseteq W_1$, $\alpha +G_1^+$
is such a cut; if necessary we also can
choose $p_1$ such that $\sign p_1=0$ ).
By \ref {obtain arbitrary signs II}
there is an elementary extension $M\succ A\< \alpha \> $ and a cut
$p$ of $M$, such that $p_1\subseteq p,\ G(p)$ lies over $G_1$ and
$V(p)$ lies over $V_1$ with $\sign p=\sign ^*\hat p=0$.
Let $W$ be a convex valuation ring of $M$, lying over $W_1$.
Then $W\subseteq V(p)$, hence again by
\ref {obtain arbitrary signs II}, there is an elementary extension
$N\succ M$ and a cut $q$ of $N$ such that
$p\subseteq q,\ G(q)$ lies over $G(p)$ and
$V(q)$ lies over $W$ with $\sign q=\sign ^*\hat q=0$.

%

\einzug
Hence $q$ is even an heir of $p$.
Moreover conditions (i), (ii) and (iv) of \ref {general main} are fulfilled
(in particular $(M,G(p))\equiv _A(N,G(q))$ by
\ref {elementary equivalent for groups}).
But $(M,p^L,V(p))\models \exists x\ p<x\in V(p)$,
$(N,q^L,V(q))\models \lnot \exists x\ q<x\in V(q)$
and $(M,p^L,V(p))\models \lnot \exists x\ V(p)<x<p$,
$(N,q^L,V(q))\models \exists x\ V(q)<x<q$.
Hence $(p,V(p))$ is not an heir of $(q,V(q))$ over $A$ and
$(q,V(q))$ is not an heir of $(p,V(p))$ over $A$.\ende


\beginlong {
\fact {Example}{}{}
Suppose $\sign p=\sign q=0$,
$Z^*(\hat p)\cap A=Z^*(\hat q)\cap A=\emptyset $,
condition (i), (ii), (iii)
and (iv)(a) of \ref {general main} holds.
Then in general $(M,p^L)\not \equiv
_A(N,q^L)$.
\pr
We construct the cuts of the following diagram from bottom to top;
again a line indicates
extensions of cuts:
$$ \diagram \hgaps {1;.8;1.5;1;1.5;1.5;1.5;2.5}\vgaps {1.3;1.3;1.3;1.3}
\object(1;1){W}
\object(1;2){}
\object(1;3){}
\object(1;4){}
\object(1;5){}
\object(1;6){}
\object(1;7){\beta \mal \Dm (W)}
\object(1;8){\beta \mal W}
\object(1;9){q}
\object(2;1){W_1}
\object(2;2){<{\beta \over \alpha }\in }
\object(2;3){V_2}
\object(2;4){}
\object(2;5){}
\object(2;6){}
\object(2;7){\beta \mal \Dm (W_1)}
\object(2;8){\beta \mal W_1}
\object(2;9){p'}
\object(3;1){}
\object(3;2){}
\object(3;3){V}
\object(3;4){{\alpha \over b}\mal V}
\object(3;5){r}
\object(3;6){\alpha \mal \Dm (V)}
\object(3;7){\alpha +r}
\object(3;8){\alpha \mal V}
\object(3;9){b\mal (\alpha +r)=:p}
\object(4;1){}
\object(4;2){}
\object(4;3){V_1}
\object(4;4){{\alpha \over b}\mal V_1}
\object(4;5){}
\object(4;6){\alpha \mal \Dm (V_1)}
\object(4;7){\alpha +{\alpha \over b}\mal V_1}
\object(4;8){\alpha \mal V_1}
\object(4;9){}
\object(5;1){}
\object(5;2){}
\object(5;3){V_0}
\object(5;4){b\ \ \ \ \ \ \ \ \ \ }
\object(5;5){}
\object(5;6){}
\object(5;7){G_0}
\object(5;8){}
\object(5;9){b\mal G_0}
\nohead\arrow1 10
\nohead\arrow7 16
\nohead\arrow8 17
\nohead\arrow16 26
\nohead\arrow9 18
\nohead\arrow12 21
\nohead\arrow17 26
\nohead\arrow18 27
\nohead\arrow21 30
\nohead\arrow22 31
\nohead\arrow23 31
\nohead\arrow24 33
\nohead\arrow25 34
\nohead\arrow34 43
\nohead\arrow26 35
\nohead\arrow27 45
\nohead\arrow30 39
\nohead\arrow33 43
\nohead\arrow35 43
\enddiagram $$\medskip

Let $V_0$ be a proper convex valuation ring of $A$ and let
$G_0$ be a convex subgroup of $A$, with $V_0\subseteq G_0$, $\sign
^*G_0^+=0$ and $V(G_0)=V_0$. Furthermore let $V_0<b\in G_0$.
Let $\alpha \models G_0^+$ and let
$V_1$ be the convex hull of $V_0$ in $A\< \alpha \> $.
Since $\sign ^*G_0^+=0$, we have that $V_1^+$ is the unique cut of
$A\< \alpha \> $ extending $V_0^+$,
$\alpha \mal \Dm (V_1)^+$ is the least and
$\alpha \mal {V_1}^+$ is the largest extension of $G_0^+$ on
$A\< \alpha \> $. Moreover $V_1^+<{\alpha \over b}\mal
V_1^+<\alpha \mal \Dm (V_1)$.\einzug
By \ref {obtain sign 0} there is some $M\succ A\< \alpha \> $
and a cut $r$ of $M$ with $\sign r=0$, such that $r$ extends
${\alpha \over b}\mal V_1^+$ and such that $G(r)$ is the convex hull
of ${\alpha \over b}\mal V_1$ in $M$.
This means that $G(r)={\alpha \over b}\mal V$, where
$V$ is the convex hull of $V_1$ in $M$.
Since $r$ extends ${\alpha \over b}\mal V_1^+$, the cut
$\alpha +r$ extends the cut $\alpha +{\alpha \over b}\mal V_1^+$.
The latter cut is $\alpha \mal (1+{1\over b}\mal V_1^+)$, so
as $b>V_1$ we have
$\alpha +{\alpha \over b}\mal V_1^+<\alpha \mal V_1^+$.
Hence
$\alpha \mal \Dm (V)<\alpha +r<\alpha \mal V^+$.
Let $p:=b\mal (\alpha +r)$. Then $G(p)=b\mal G(\alpha +r)=
b\mal G(r)=\alpha \mal V$ and $\sign p=0$.\einzug
Let $\beta $ be a realization of $\alpha \mal V^+$.
Since $\sign p=0$, the unique cut $p'$ of $M\< \beta \> $ extending
$p$ has again signature 0 and $G(p')$ is the largest convex subgroup
of $M\< \beta \> $ lying over $G(p)=\alpha \mal V$.
Hence $G(p')=\beta \mal V_2$, where $V_2$ is the largest convex
valuation ring of $M\< \beta \> $ lying over $V$.
let $W_1\proper V_2$ be a convex valuation ring of $M\< \beta \> $
(for example we can take
$W_1$ to be the convex hull of $V$ in $M\< \beta \> $).
In particular $\beta \mal W_1\subseteq G(p')$.\einzug
Finally let $N\succ M$ and let $q$ be a cut of $N$ with $\sign q=0$ such
that $G(q)$ is the convex hull of $\beta \mal W_1$ in $N$
(by \ref {obtain sign 0}). Hence $G(q)=\beta \mal W$, where
$W$ is the convex hull of $W_1$ in $N$.\einzug
\claim $\beta \mal \Dm (W_1)^+$ extends $\alpha \mal V^+$.

It is enough to prove the claim if $W_1$ is the convex hull of $V$.
Since
$W_1<\eps :={\beta \over \alpha } \in V_2$, we get
${1\over \eps }\mal \beta \mal W_1^+=
\alpha \mal W_1^+$ extends $\alpha \mal V^+$ and
${1\over \eps }\mal \beta \mal W_1^+<\beta \mal \Dm (W_1)^+$.\\

The claim implies that $\beta \mal \Dm (W)^+$ extends $\alpha \mal V^+$.
In particular ${1\over b}\mal q<\beta \mal \Dm (W)^+$ and
${1\over b}\mal q\mal V(q)^+<\hat q$.
On the other hand ${1\over b}\mal p\mal V(p)^+=\hat p$, which proves
$$(N,q^L)\not \equiv _A(M,p^L)$$

But $(N,G(q))\equiv _A(M,G(p))$
($\sign _A^*\alpha V=\sign _A^*\beta W=0$),
$q$ is an heir of $p$
and for all $A$-definable maps
$f$ we have $f(p)<V(p)^+$ iff $f(q)<V(q)^+$ ($\sign ^*G_0^+=0$).
Note also that we can choose $W_1$ to be the convex hull of $V$, hence
$V(q)=W$ lies over $V(p)=V$, $\hat q$ lies over $\hat p$ and $q$ is an
heir of $p$.

\medskip
}\endlong

\beginlong {
\fact {example}{need condition (d)}{}
The following example shows that we need condition
(d) in \ref {general main}.\einzug
We start with a model $A$ of $T$, a convex subgroup
$G_0$ of $A$ with $\sign ^*G_0^+=0$ and an element
$b\in A$ with $G_0^+<b\in \Dm (G_0)$.
Let $\beta $ be a realization of $G_0^+$,
$M_1:=A\< \beta \> $, $V_1:=$ the convex hull of $V(G_0)$ in
$M_1$ and let $G_1:=\beta \mal V_1$.
Let $M_2\succ M_1$ such that there are realizations $\alpha _2,\beta _2$
of $V_1^+$ together with a convex valuation ring $V_2$ of $M_2$
with $\alpha _2\in V_2<\beta _2$.
Clearly

\underconstruction{noch unklar}
}\endlong

\beginlong {
\OP {}{
Here are some more problems.
\beg {(a)}
\lsmall Suppose in \ref {qe language}
we do not put all the $S_\lambda ^k$ and
$R_\eta ^k$ into the new language.
Prove that we don't have qe for some $\overline {T_\delta ^\eps}$.\\
(this would probably also imply that qe can not be achieved by
adjoining finitely many new symbols).
\lsmall
Prove the box theorem with parameters (see also (d) below).
\lsmall Say $\sign p=\sign ^*\hat p=0$ here.
Suppose we know (i),(ii) and (iv) of \ref {all signs 0} and (iii) for
$\lambda =1$. Is this enough to get $(M,p^L)\equiv _A(N,q^L)$ ?\herz 1\\
Observe that $$\{ \beta \in M\st a+\beta \mal \hat p<p \} ^+={1\over
\widetilde {p-a}}.$$
\lsmall
Let $G\subseteq M$ be a convex subgroup and let
$a_1,...,a_n\in M$.
Describe the subsets of $M$, which are 0-definable in
$(M,a_1G,...,a_nG)$.
Which $\{ a_1,...,a_n\} $-definable subsets of
$(M,G)$ are 0-definable in $(M,G,a_1G,...,a_nG)$ ?
\lsmall Prove $(M,p^L)\equiv _A(N,q^L)\iff (M,p^L)\equiv _{\{ a,b\} }
(N,q^L)$ for all $a,b\in A$ in the case where
$V(p)$ and $V(q)$ are $T$-convex and $T$ is exponentially
bounded with levels, say for $T_{an,exp}$.
Prove this at least for edges  of convex subgroups.
\endbeg
}

\subtitle{Sets definable in $(M,p^l)$ with parameters from $M$.}

Here we may assume that $\sign p=0$;

Let $M\prec A\prec N$ and let $q$ be a cut of $M$
such that $(M,p^L)\prec (N,q^L)$.
Let $p_0:=q\uhr A$, $G_0^+:=\hat q\uhr A$ and $V_0:=V(q)\cap A$.\\
\claim 1 Suppose $p_0=a+\hat p_0$ with $\alpha \in A$.
}\endlong

\bigskip
\zl{{\sc References}}\smallskip
\beg{---------------}
\lit{vdD1}{L. van den Dries}{Tame Topology and O-minimal Structures}{
London Mathematical Society Lecture Note Series, 248}
\lit{vdD2}{L. van den Dries}{$T$-convexity
and tame extensions II}{J. Symbolic Logic 62 (1997),no. 1, 14-34}
\lit{vdD-L}{L. van den Dries, A.H. Lewenberg}{$T$-convexity
and tame extensions}{Journal of Symbolic Logic, vol. 60, no. 1, 74-101 (1995)}
\lit{vdD-S}{L. van den Dries, P. Speissegger}{The field of reals with
multisummable series and the exponential function}{Preprint 1997}
\lit{MMS}{D. Macpherson, D. Marker, C. Steinhorn}{Weakly o-minimal
structures and real closed fields}{ Trans. Amer. Math. Soc. 352 (2000), no. 12, 5435--5483}
\lit {Ma}{D. Marker}{Omitting types in $o$-minimal theories}
{J. Symb. Logic, vol.51, no 1, 1986}
\lit{Mi1}{C. Miller}{A growth dichotomy for o-minimal expansions of ordered
fields}{Logic: from foundations to applications (European Logic Colloquium,
1993, W. Hodges, J. Hyland, C. Steinhorn and J. Truss, editors),
Oxford University Press, 1996, pp. 385-399.}
\lit{Mi2}{C. Miller}{Exponentiation is hard to avoid}
{Proc. Am. Math. Soc., vol. 122, no. 1, 257-259 (1994)}
\lit {Poi}{Poizat, B.}{Cours de th\'eorie des mod\`eles}
{Villeurbanne: Nur al-Mantiq wal-Ma'rifah (1985)}
\lit{Tr}{M. Tressl}{Model Completeness of o-minimal Structures
expanded by Dedekind Cuts}{submitted; electronically available under
http://www-nw.uni-regensburg\\
.de/~.trm22116.mathematik.uni-regensburg.de/papers/cutsa.ps}
\lit{Wi}{A.J. Wilkie}{Model completeness results for expansions of the
ordered field of real numbers by restricted pfaffian functions and the
exponential function}{J. Am. Math. Soc., vol 9, Number 4, 1051-1094 (1996)}
\endbeg\medskip

\bye